\documentclass[11pt]{amsart}
\usepackage{amsmath,amsfonts,latexsym,amssymb,amsthm,mathrsfs,array}

\usepackage{amsfonts}
\usepackage[usenames]{color}
\usepackage[colorlinks=true,linkcolor=webgreen,filecolor=webbrown,citecolor=webgreen]
{hyperref}
\definecolor{webgreen}{rgb}{0,.5,0}
\definecolor{webbrown}{rgb}{.6,0,0}

\newcommand{\ZZ}{{\mathbb Z}}

\def\ZZ{{\mathbb Z}}

\newtheorem{dfn}{Definition}[section]
\newcommand{\bdfn}{\begin{dfn}\rm}
	\newcommand{\edfn}{\end{dfn}}
\newtheorem{thm}[dfn]{Theorem}
\newcommand{\bthm}{\begin{thm}}
	\newcommand{\ethm}{\end{thm}}
\newtheorem{lmma}[dfn]{Lemma}                   
\newcommand{\blmma}{\begin{lmma}}                   
	\newcommand{\elmma}{\end{lmma}}                   
\newtheorem{ppsn}[dfn]{Proposition}
\newcommand{\bppsn}{\begin{ppsn}}
	\newcommand{\eppsn}{\end{ppsn}}
\newtheorem{crlre}[dfn]{Corollary}
\newcommand{\bcrlre}{\begin{crlre}} 
	\newcommand{\ecrlre}{\end{crlre}}
\newtheorem{rmk}[dfn]{Remark}
\newcommand{\brmk}{\begin{rmk}\rm} 
	\newcommand{\ermk}{\end{rmk}}

\numberwithin{equation}{section}

\title[Hamiltonian Extended Affine Lie Algebra]{Hamiltonian Extended Affine Lie Algebra and its Representation Theory}
\author{S. Eswara Rao}
\address{S. Eswara Rao, School of Mathematics, Tata Institute of Fundamental Research, Homi Bhaba Road, Colaba, Mumbai 400005, India.}
\email{senapati@math.tifr.res.in, sena98672@gmail.com}
\date{}

\begin{document}
	
\maketitle
\begin{abstract}
We introduce a new class of extended affine Lie algebras called Hamiltonian Extended Affine Lie Algebras (HEALAs). They are so called because the corresponding derivation algebra is the classical Hamiltonian algebra. We classify the irreducible integrable modules for HEALA based on the classification of irreducible Jet modules for the Hamiltonian algebra (both zero level and non-zero level).\\\\
{\bf{2020 MSC}:} Primary: 17B65; Secondary: 17B67, 17B70.\\
{\bf{KEY WORDS}:} Hamiltonian algebra, integrable modules, extended affine Lie algebras, triangular decomposition, jet modules, toroidal Lie algebras.
		
\end{abstract}
	
\section{Introduction}
In this paper, we introduce a new class of extended affine Lie algebras (EALAs) called the Hamiltonian extended affine Lie algebras (HEALAs) and study their representation theory. EALAs have been extensively studied in the last two decades (see \cite{AABGP,AG,ABFP,N1,N2,Y} and the references therein). These EALAs are generalizations of finite-dimensional simple Lie algebras (nullity $0$) and affine Lie algebras (nullity $1$). The representation theory of both finite-dimensional Lie algebras and affine Lie algebras are extremely well-developed classical objects. The representation theory of toroidal Lie algebras are still in progress \cite{MEY,ES,EM,E,EJ,ESB}. Most often, integrable modules are classified for the full toroidal Lie algebra and its related subalgebras. The most important EALA is the toroidal extended affine Lie algebra (TEALA) which is reasonably well-developed. Integrable modules are very important subclass of representations. The irreducible integrable modules over TEALA of non-zero level and nullity at least $3$ were classified in \cite{ESB}. The nullity $2$ case was settled in \cite{CLT1,CLT2} only for the non-twisted case.  	

In this paper, we study irreducible integrable modules for HEALA where the core of the HEALA acts non-trivially (see Section \ref{S3} for the definition of HEALA). We use a very important and powerful classification theorem for Jet modules over the Hamiltonian Lie algebra given in \cite{T}. If the core acts trivially, the classification problem gets reduced to the classification of irreducible modules for the Hamiltonian Lie algebra which is still open.

We now give details of each section. In Section \ref{S2}, we recall the toroidal Lie algebra, its root system and the notion of integrable modules. The root systems of toroidal Lie algebras, full toroidal Lie algebras, TEALAs and HEALAs are all the same. 

In Section \ref{S3}, we define the HEALA as a subquotient of TEALA. The corresponding derivation algebra is a Hamiltonian algebra and hence the name Hamiltonian EALA. We also define the twisted version of HEALA but the representation theory is postponed for the next publication.

In Section \ref{S4}, we consider the two variable case and classify the irreducible integrable modules for HEALA of non-zero level (Theorem \ref{T4.4}). The level zero case is presented in Section \ref{S7}. We note that the case of two variables is also studied in \cite{CLT1,CLT2}, but their triangular decomposition is different from ours. It will be interesting to compare both these classes of modules. The techniques we use in the two variable case will be very useful in Section \ref{S5} and Section \ref{S6}. In Section \ref{S8}, we define yet another EALA and call it KEALA and this time the number of variables is odd. 	

\section{Toroidal Lie algebras and full toroidal Lie algebras} \label{S2}
In this section, we recall the definitions of toroidal Lie algebras and full toroidal Lie algebras.\\
(2.1) We fix a positive integer $N$ throughout this paper. All the vector spaces, algebras and tensor products are over the field of complex numbers $\mathbb{C}$. We shall denote the set of integers, natural numbers and non-negative numbers by $\mathbb{Z}$, $\mathbb{N}$ and $\mathbb{Z}_{+}$ respectively. Let $\mathbb{Z}^N = \bigoplus \mathbb{Z}$ stand for $N$-copies of $\mathbb{Z}$, $r = (r_1, \ldots, r_N) \in \mathbb{Z}^N$ and similarly $s,k,l$ denote vectors in $\mathbb{Z}^N$. Let $A = A_N = \mathbb{C}[t_1^{\pm 1}, \cdots, t_N^{\pm 1}]$  be the algebra of Laurent polynomials in $N$ commuting variables. For $r \in \mathbb{Z}^N$, let $t^{r}=t_1^{r_1} \cdots t_N^{r_N} \in A$.\\
(2.2) Let $\mathfrak{g}$ be a finite-dimensional simple Lie algebra. Let $\mathfrak{h}$ be a Cartan subalgebra of $\mathfrak{g}$ and $\mathfrak{g} = \bigoplus_{\alpha \in \Delta} \mathfrak{g}_{\alpha} \oplus \mathfrak{h}$  
be the root space decomposition of $\mathfrak{g}$ relative to $\mathfrak{h}$, where $\Delta$ is the corresponding (finite) root system.
Let $<\cdot,\cdot>$ be a non-degenerate, symmetric and invariant bilinear form on $\mathfrak{g}$.

Let $\Omega_A$ denote the vector space spanned by the symbols of the form $\{t^rK_i  \ | \ r \in \mathbb{Z}^N, \ 1 \leqslant i \leqslant N \}$. It is clear that $\Omega_A$ is $\mathbb{Z}^N$-graded, with each component being $N$-dimensional. Let $d_A$ be the vector subspace spanned by  $\{\sum_{i =1}^{N}r_i t^rK_i, \ r \in \ZZ^{N} \}$. Then let $\mathcal{Z} = \Omega_A / d_A$ which is $\mathbb{Z}^N$-graded and each component is $(N-1)$-dimensional, except the zeroth component that is $N$-dimensional. Note that in the particular case $N=1$, $\mathcal{Z}$ is $1$-dimensional. Let $\mathbb{C}^N = \bigoplus \mathbb{C}$ stand for $N$-copies of $\mathbb{C}$ and take $u = (u_1, \ldots, u_N) \in \mathbb{C}^N$. There is a standard non-degenerate bilinear form on $\mathbb{C}^N$, which we shall denote by $(\cdot, \cdot)$. Let $K(u,r) = \sum_{i=1}^{N}u_i t^rK_i, \ u \in \mathbb{C}^N, \ r \in \mathbb{Z}^N$. It is clear that $K(r,r) \in d_A$.\\ 
(2.3) We shall now define toroidal Lie algebra given by
	\begin{align*}
		\tau = \mathfrak{g} \otimes A \oplus \mathcal{Z} \oplus D, \ \text{where}
	\end{align*} 
$D$ is spanned by $\{d_1, \ldots, d_N \}$ which are the degree derivations.\\
Let $X(r) = X \otimes t^r, \ X \in \mathfrak{g}, \ r \in \mathbb{Z}^N \ \text{and} \ \mathfrak{g}(r) = \mathfrak{g} \otimes \mathbb{C}t^r$. The Lie bracket is given by \\
(2.3)(1) $[X(r), Y(s)] = [X,Y](r + s) + <X,Y> K(r,r+s)$.\\
(2.3)(2) $\mathcal{Z}$ is central in $\mathfrak{g} \otimes A$.\\
(2.3)(3) $[d_i, X(r)] = r_i X(r)$, $[d_i, K(u,r)]= r_iK(u,r)$, $[d_i, d_j] = 0$\\
for all $r,s \in \mathbb{Z}^N, \ u \in \mathbb{C}^N, \ X,Y \in \mathfrak{g}, \ 1 \leqslant i,j \leqslant N$.\\
(2.4) Let $\widetilde{\mathfrak{h}} = \mathfrak{h} \bigoplus \mathcal{Z}_0 \bigoplus D$, where we have \\
 $\mathcal{Z}_0 = \text{span} \{K(u,0) \ | \ u \in \mathbb{C}^N \} = \text{span} \{K_1, \ldots, K_N \}$. In this case, $\widetilde{\mathfrak{h}}$ will serve as a Cartan subalgebra of $\tau$.\\
(2.5) Let $DerA$ be the derivation algebra of $A$. It is well-known that\\
$DerA = \text{span} \{t^rd_i \ | \ r \in \mathbb{Z}^N, \ 1 \leqslant i \leqslant N \}$, where $d_i = t_i \dfrac{\partial}{\partial t_i}$. Let \\
$D(u,r) = \sum_{i=1}^{N} u_it^rd_i, \ u \in \mathbb{C}^N, \ r \in \mathbb{Z}^N$. Then we have \\
$[D(u,r), D(v,s)] = D(w, r+s)$, where $w = (u,s)v - (v,r)u$. It is known that $DerA$ admits an abelian extension on $\mathcal{Z}$ (see \cite{EM}) as \\
(2.5)(1) $[D(u,r), D(v,s)] = D(w, r+s) + (u,s)(v,r) K(r,r+s)$.\\
(2.5)(2) $[D(u,r), K(v,s)] = (u,s)K(v,r+s) + (u,v)K(r,r+s)$.\\ 
It is known that $DerA$ has no non-trivial central extension for $N \geqslant 2$ \cite{RSS}. For $N=1$, the above abelian extension becomes a central extension.\\
(2.6) We shall now define the full toroidal Lie algebra given by
\begin{align*}
\widetilde{\tau} = \mathfrak{g} \otimes A \oplus \mathcal{Z} \oplus DerA.
\end{align*}
The Lie bracket is given by (2.3)(1), (2.3)(2), (2.5)(1), (2.5)(2) and the following.\\
(2.6)(1) $[D(u,r), X(s)] = (u,s)X(r+s)$.\\
(2.7) We shall now define automorphisms of $\widetilde{\tau}$. Let $GL(N, \mathbb{Z})$ be the group of integral matrices with determinant $\pm 1$. Let $B \in GL(N, \mathbb{Z})$ which naturally acts on $\mathbb{Z}^N$ and we denote this action by $Bu, \ u \in \mathbb{C}^N$. We shall now define an automorphism of $\widetilde{\tau}$, again denoted by $B$, by setting
\begin{align*}
	B.X(r) = X(Br), \ B.K(u,r) = K(Bu,Br), \ B.D(u,r) = B(Fu,Br), 
\end{align*}
where $F = (B^T)^{-1}$ (see \cite{EJ} for details). It can be directly verified that $B$ defines an automorphism of $\widetilde{\tau}$.\\
(2.8) We shall now define roots and co-roots of $\widetilde{\tau}$ and the corresponding Weyl group. We shall also define integrable modules for $\widetilde{\tau}$ which is the main topic of this paper. As the root system of $\widetilde{\tau}$ is same as that of the toroidal Lie algebra $\tau$, we simply recall the notions from \cite{E}.

Let $\delta_i$ ($1 \leqslant i \leqslant N$) be in $\widetilde{\mathfrak{h}}^*$ defined by $\delta_i(h) = 0, \ \delta_i(\mathcal{Z}_0) = 0$ and $\delta_i(d_j) = \delta_{ij}$. Also denote $\omega_i$ ($1 \leqslant i \leqslant N$) in $\mathfrak{\widetilde{h}}^*$ by $\omega_i(h) = 0, \ \omega_i(d_j) = 0$ and $\omega_i(K_j) = \delta_{ij}$. For $r \in \mathbb{Z}^N$, let $\delta_r = \sum_{i=1}^{N} r_i \delta_i$. Let $\{\alpha_1, \ldots, \alpha_d \}$ be a set of simple roots of $\mathfrak{g}$ where $d = \text{dim} \ \mathfrak{h}$. Let $\alpha_1^{\vee}, \ldots, \alpha_d^{\vee}$ be the corresponding co-roots. Then $\{\alpha_1, \ldots, \alpha_d, \delta_1, \ldots, \delta_N, \omega_1, \ldots, \omega_N \}$ is a basis of $\widetilde{\mathfrak{h}}^*$. Now $\widetilde{\mathfrak{h}}^*$ admits a symmetric, non-degenerate bilinear form such that $<\alpha, \delta_i> = 0 = <\alpha, \omega_i> \ \forall \ \alpha \in \Delta, \ <\delta_i, \delta_j> = 0 = <\omega_i, \omega_j>, \ <\delta_i, \omega_j> = \delta_{ij} \ \forall \ 1 \leqslant i,j \leqslant N$. Again $\widetilde{\mathfrak{h}}$ has a basis $\{\alpha_1^{\vee}, \ldots, \alpha_d^{\vee}, K_1, \ldots, K_n, d_1, \ldots, d_N \}$ with the bilinear form $<h,d_i> = 0 = <h,K_i>, \ <d_i,d_j> = 0 = <K_i,K_j>, \ <d_i, K_j> = \delta_{ij}$. On $\mathfrak{h}$, the form is simply the restriction of $<\cdot, \cdot>$ onto $\mathfrak{h}$. We certainly have dim $\widetilde{\mathfrak{h}} = \widetilde{\mathfrak{h}}^* = d + 2N$.\\
(2.9) Let $\Delta_{re} = \{\alpha + \delta_r \ | \ \alpha \in \Delta, \ r \in \mathbb{Z}^N \}$ and $\Delta_{im} = \{\delta_r \ | \ r \in \mathbb{Z}^N \}$. Then $\widetilde{\Delta} = \Delta_{re} \cup \Delta_{im}$ is the root system of $\widetilde{\tau}$ with respect to the Cartan subalgebra $\widetilde{\mathfrak{h}}$. Note that $(\delta_r, \delta_s) = 0$ and $(\alpha + \delta_r, \beta + \delta_s) = (\alpha, \beta) \ \forall \ \alpha, \beta \in \Delta$.\\
For $\lambda \in \widetilde{\mathfrak{h}}^*$, let $\overline{\lambda}$ denote the restriction of $\lambda$ to $\mathfrak{h}$. And for any $\mu \in \mathfrak{h}^*$, we can extend to $\widetilde{\mathfrak{h}}^*$ by simply setting $\mu(d_i) = 0 = \mu(K_i)$ so that any $\lambda \in \widetilde{\mathfrak{h}}^*$ can be expressed as $\lambda = \overline{\lambda} + \sum_{i=1}^{N}g_i \delta_i + \sum_{i=1}^{N} s_i \omega_i$. We also define $\alpha_{d + j} = -\beta + \delta_j, \ 1 \leqslant j \leqslant N$, where $\beta$ is the highest root of $\mathfrak{g}$. Then the collection $\{\alpha_1, \ldots, \alpha_d, \alpha_{d+1}, \ldots, \alpha_{d+N}, \omega_1, \ldots, \omega_N \}$ also forms a basis of $\widetilde{\mathfrak{h}}^*$. One can define co-root $\gamma^{\vee}$ for $\gamma \in \Delta_{re}$ (see \cite{E}).\\
(2.10) We define an ordering on $\widetilde{\mathfrak{h}}^*$. For $\lambda, \mu \in \widetilde{\mathfrak{h}}^*$, we call $\mu \leqslant \lambda$ if $\lambda - \mu = \sum_{i=1}^{d} n_i \alpha_i + n_{d + m} \delta_m + n_{d + 2m} \delta_{2m}$, where all the $n_i$'s are integers.\\
Either (i) $n_{d+m} - n_{d+2m} > 0$ or (ii) $n_{d+m} = n_{d+2m} > 0$ or (iii) $n_{d+m} = n_{d+2m} = 0$ but $\sum_{i=1}^{d} n_i \alpha_i \in Q^+$, where $Q^+$ is the span of all non-negative linear combinations of $\alpha_1, \ldots, \alpha_d$.\\
(2.11) For any real root $\gamma$, define the reflection operator $r_{\gamma}$ on
$\widetilde{\mathfrak{h}}^*$ by setting
\begin{align*}
r_{\gamma}(\lambda) = \lambda - \lambda(\gamma^{\vee}) \gamma.
\end{align*}    
Then we define the Weyl group $W$ as the group generated by the reflections $r_{\gamma}, \ \gamma \in \Delta_{re}$ (see \cite{E}).\\
(2.12) A module $V$ over $\widetilde{\tau}$ (or over $\tau$) is called integrable if
\begin{enumerate}
	\item $V = \bigoplus_{\lambda \in \widetilde{\mathfrak{h}}^*} V_{\lambda}$, where $V_{\lambda} = \{v \in V \ | \ hv = \lambda(h)v \ \forall \ h \in \widetilde{\mathfrak{h}}^* \}$. 
	\item $X_{\alpha}(r)$ acts locally nilpotently on $V$ for all $\alpha \in \Delta, \ r \in \mathbb{Z}^N, \ X_{\alpha} \in \mathfrak{g}_{\alpha}$. 
\end{enumerate}
(2.13) Let $V$ be an integrable module for $\widetilde{\tau}$ (or for $\tau$) with finite-dimensional weight spaces (i.e. dim $V_{\lambda} < \infty \ \forall \ \lambda \in \widetilde{\mathfrak{h}}^*$). We shall denote the set of all weights of $V$ by $P(V) = \{ \lambda \in \widetilde{\mathfrak{h}}^* \ | \ V_{\lambda} \neq 0 \}$. Then the following are very standard (see \cite{E}).\\
(2.13)(1) $P(V)$ is $W$-invariant.\\
(2.13)(2) $\text{dim} V_{\lambda} = V_{w \lambda} \ \forall \ w \in W, \ \lambda \in P(V)$.\\
(2.13)(3) If $\alpha \in \Delta_{re}$ and $\lambda \in P(V)$, then $\lambda(\alpha^{\vee}) \in \mathbb{Z}$.\\
(2.13)(4) If $\alpha \in \Delta_{re}, \lambda \in P(V)$ and $\lambda(\alpha^{\vee}) > 0$, then $\lambda - \alpha \in P(V)$.\\
(2.13)(5) $\lambda(K_i)$ is a constant integer for all $\lambda \in P(V)$, if $V$ is irreducible.\\
(2.14) We use the following terminology in the later sections.\\
Let $\tau$ be any Lie algebra and suppose that $\tau = \tau_{-} \oplus \tau_0 \oplus \tau_+$ is called a triangular decomposition if $\tau_+, \ \tau_{-}$ and $\tau_0$ are subalgebras and $[\tau_0, \tau_{\pm}] \subseteq \tau_{\pm}$. Let $W$ be any irreducible module for $\tau_0$. Then consider the Verma module 
\begin{align*}
M(W) = U(\tau) \bigotimes_{\tau_+ \oplus \tau_0} W,
\end{align*}  
 where $\tau_+$ acts trivially on $W$. Let $L(W)$ be the unique irreducible quotient of $M(W)$ (if it exists, all our triangular decompositions have this property). In general, we call $W$ as the top of $M(W)$ and $L(W)$ is the corresponding irreducible quotient. It is clear that $L(W_1) \cong L(W_2)$ as $\tau$-modules if and only if $W_1 \cong W_2$ as $\tau_0$-modules.\\
 (2.15) We say that a module has level zero if the zero degree central elements $K_i$ act trivially for all $1 \leqslant i \leqslant N$, otherwise it is said to have non-zero level.  	
\section{Hamilton Extended Affine Lie Algebra (HEALA)} \label{S3}
In this section, we introduce a new class of extended affine Lie algebras called HEALAs. We shall first introduce the Hamiltonian Lie algebra which is a classical object and a subalgebra of $DerA$. It is also called $H$ type.\\
(3.1) We shall closely follow \cite{T} to define the Hamiltonian Lie algebra and denote it by $H_N$. We need a bit of preparation for this. The Hamiltonian Lie algebras are known to exist only when $N$ is even. Thus we consider $N = 2m$, with $m \geqslant 1$.\\
Let $r = (r_1, r_2, \ldots, r_m, r_{m+1}, \ldots, r_{2m}) \in \mathbb{Z}^N$ and also define \\
$\overline{r} = (r_{m+1}, \ldots, r_{2m}, -r_1, -r_2, \ldots, -r_m)$. It is easy to check the following.\\
(3.1)(1) $(\overline{r},r) = 0$.\\
(3.1)(2) $(r, \overline{s}) = - (\overline{r},s), \ \overline{\overline{r}} = -r \ \forall \ r,s \in \mathbb{Z}^N$.\\
(3.1)(3) For $r \in \mathbb{Z}^N$, let us set\\
$h_r = D(\overline{r},r) = \sum_{i=1}^{m} (r_{m+i}t^rd_i - r_it^rd_{m+i})$ (see (2.5)).\\
One can check that $[h_r,h_s] = (\overline{r},s)h_{r+s}, \ [h_r, h_{-r}] = 0, \ h_0 = 0$.\\
As earlier, $[D(u,0), h_r] = (u,r)h_r$. Let $H_N = \text{span} \{h_r \ | \ 0 \neq r \in \mathbb{Z}^N \}$ and $\widetilde{H_N} = H_N \oplus D$. Then $H_N$ and $\widetilde{H_N}$ are known to be simple.\\
(3.2) We now recall the definition of an extended affine Lie algebra (EALA).\\
Let L be any Lie algebra.\\
(EA1) L possesses a non-degenerate symmetric bilinear form $(\cdot | \cdot)$ which is invariant in the sense that $([x,y] | z) = (x | [y,z]) \ \forall \ x,y,z \in \text{L}$.\\
(EA2) L has a non-trivial finite-dimensional self-centralizing ad-diagonalizable abelian subalgebra H.\\
We shall assume three further axioms about the triplet (L, $(\cdot | \cdot)$, H). To describe these axioms, we need some more notations. Using (EA2), we have 
\begin{align*}
\text{L} = \bigoplus_{\alpha \in H^*} \text{L}_{\alpha}, \ \text{where} \ \text{L}_{\alpha} = \{x \in \text{L} \ | \ [h,x]= \alpha(h)x \ \forall \ h \in \text{H} \}.	
\end{align*}
Let $R = \{\alpha \in \text{H}^* \ | \ \text{L}_{\alpha} \neq 0 \}$ which is called the root system of L with respect to H. We have $0 \in R$. Also $\alpha, \beta \in R$ and $\alpha + \beta \neq 0 \implies (\text{L}_{\alpha}, \text{L}_{\beta}) = 0$.\\
Thus $R = -R$. The form restricted to H is non-degenerate and so we can transfer the form to H$^*$. Let 
\begin{align*}
R^{\times} = \{\alpha \in R \ | \ (\alpha | \alpha) \neq 0 \}, \ R^{0} = \{\alpha \in R \ | \ (\alpha | \alpha) = 0 \}.  
\end{align*}
The elements of $R^{\times}$ (respectively $R^0$) are called non-isotropic (respectively isotropic) roots. We also have $R = R^{\times} \cup R^0$.\\
(EA3) If $x_{\alpha} \in \text{L}_{\alpha}, \ \alpha \in R^{\times}$, then ad$x_{\alpha}$ should act locally nilpotently on L.\\
(EA4) $R$ is a discrete subset of H$^*$.\\
(EA5) $R$ is an irreducible root system, which means that\\
(a) $R^{\times} = R_1 \cup R_2$ and $(R_1, R_2) = 0$ should imply either $R_1 = 0$ or $R_2 = 0$.\\
(b) If $\sigma \in R^0$, then there exists $\alpha \in R^{\times}$ such that $\sigma + \alpha \in R$.\\
Any Lie algebra satisfying (EA1) to (EA5) is called an EALA.\\
Extensive research has been done on the structure and classification of EALAs. See \cite{AABGP,AG,ABFP,N1,N2} and also the references therein.\\
(3.3) We shall now give some well-known examples of EALAs. The full toroidal Lie algebra $\widetilde{\tau}$ falls short of an EALA, as it does not satisfy (EA1).\\
(3.4) A certain subalgebra of $\widetilde{\tau}$ is an EALA. Define
\begin{align*}
S_N = \text{span} \{D(u,r) \ | \ (u,r) = 0, \ u \in \mathbb{C}^N, \ r \in \mathbb{Z}^N \},
\end{align*}
which is a subalgebra of $DerA$. It is known to be simple and one of the four infinite sereis. It is called $S$ type. Note that $S_2 = \widetilde{H_2}$. Let
\begin{align*}
\tau(S_N) = \mathfrak{g} \otimes A \oplus \mathcal{Z} \oplus S_N.
\end{align*}
Define a bilinear form $(\cdot | \cdot)$ on $\tau(S_N)$ as follows.\\
(3.4)(1) $(X(r) | Y(s)) = <X,Y> \delta_{r+s,0}$, \\
(3.4)(2) $(D(u,r) | K(v,s)) = \delta_{r+s,0}(u,v) \ \forall \ X, Y \in \mathfrak{g}, \ r,s \in \mathbb{Z}^N, \ u,v \in \mathbb{C}^N$.\\
All other values are zero. It is a standard fact that $\tau(S_N)$ is an EALA with the above non-degenerate bilinear form.\\
(3.5) We shall now define HEALA. We need some preparation for that.\\
Let $\mathcal{K} = \text{span} \{K(u,r) \in \mathcal{Z} \ | \ r \neq 0, \ (u, \overline{r}) = 0 \}$. Clearly $\mathcal{K}$ is $\mathbb{Z}^N$-graded, with each of its graded components being $(N-2)$-dimensional. Observe that $\mathcal{K}$ is trivial if $N = 2$.\\
(3.6) It is straightforward to check that $[H_N, \mathcal{K}] \subseteq \mathcal{K}$ using (2.5)(2). Now consider the Lie algebra given by
\begin{align*}
\tau^{\prime} = \mathfrak{g} \otimes A \oplus \mathcal{Z} \oplus H_N \oplus D.
\end{align*}
Note that $\mathcal{K}$ is an ideal in $\tau^{\prime}$ and thus
\begin{align*}
\tau(H_N) = \mathfrak{g} \otimes A \oplus \mathcal{Z}/\mathcal{K} \oplus H_N \oplus D
\end{align*}
is also a Lie algebra. Clearly $\mathcal{Z}/\mathcal{K}$ is $\mathbb{Z}^N$-graded and each component is $1$-dimensional except the zeroth component which is $N$-dimensional.

\bppsn
\
\begin{enumerate}
\item $\tau(H_N)$ is an EALA.
\item $\text{dim} (\mathcal{Z}/\mathcal{K})_r = 1 \ \forall \ r \neq 0$.
\item $\text{dim} (\mathcal{Z}/\mathcal{K})_0 = N$. 
\end{enumerate}
\eppsn

\begin{proof}
Suppose for $r \neq 0, \ (\overline{r}, u) \neq 0$ and $(\overline{r}, v) \neq 0$ for some $u, v \in \mathbb{C}^N$. Then
\begin{align*}
\bigg (\dfrac{u}{(\overline{r},u)} - \dfrac{v}{(\overline{r}, v)} \ , \ \overline{r} \bigg) = 0. 
\end{align*}
Thus $K(u,r) = \lambda K(v,r) \in \mathcal{Z}/ \mathcal{K}$ for some $\lambda \neq 0$. This proves (2). We can take $K(\overline{r}, r)$ as a basis of $\mathcal{Z}/ \mathcal{K}$ along with $\{K_1, \ldots, K_N \}$. \\
To prove (1), it is easy to verify all the axioms of EALA are true except the bilinear form (EA1). Define
\begin{align*}
(X(r) \ | \ Y(s)) = <X,Y> \delta_{r+s,0} \ \forall \ r,s \in \mathbb{Z}^N, \\ 
(D(\overline{r}, r) \ | \ K(\overline{s}, s) ) = \delta_{r+s,0}(\overline{r}, \overline{s}) \ \forall \ 0 \neq r,s \in \mathbb{Z}^N, \\
(D(u,0) \ | \ K(v,0)) = (u,v) \ \forall \ u,v \in \mathbb{C}^N. 
\end{align*}
All the other values are zero. \\
Note that $(D(\overline{r}, r) | K) = 0$. We need to check the following.\\
$([D(\overline{r}, r), D(\overline{s}, s)] \ | \ K(\overline{l}, l) ) = (D(\overline{r}, r) \ | \ [D(\overline{s}, s), K(\overline{l}, l)])$ and \\
$([D(\overline{r}, r), D(\overline{s}, s)] \ | \ D(\overline{l}, l) ) = (D(\overline{r}, r) \ | \ [D(\overline{s}, s), D(\overline{l}, l)])$ if $l+r+s =0$.\\
But both of them follow from direct checking.\\
Now to check that the form is non-degenerate, consider for $r+s = 0$,\\
$(D(\overline{r}, r) | K(\overline{s}, s) ) = \delta_{r+s,0}(\overline{r}, \overline{s}) = -(\overline{s}, \overline{s}) = - (s,s) \neq 0$. It is clear that the bilinear form is actually descending from the form on $\tau(S_N)$ and we already know that $(D(\overline{r}, r) \ | \ K) = 0$.
\end{proof}

\noindent (3.7) Recall that $B \in GL(N, \mathbb{Z})$ (see (2.7)) defines an automorphism of $\widetilde{\tau}$. In fact, it leaves $\tau(S_N)$ invariant. Now define\\
$\mathcal{K}_B = \text{span} \{K(Br,Bs) \ | \ (r, \overline{s}) = 0 , \ r,s \neq 0 \}$ and \\
$H_B = \text{span} \{D(E \overline{r},Bs) \ | \ 0 \neq r \in \mathbb{Z}^N \}$, where $E = (B^T)^{-1}$.\\
One can check that $[H_B, \mathcal{K}_B] \subseteq \mathcal{K}_B$. Consider 
\begin{align*}
\tau_B = \mathfrak{g} \otimes A \oplus \mathcal{Z} \oplus H_B \oplus D.	
\end{align*}
Then $\mathcal{K}_B$ is an ideal of $\tau_B$. This implies that
\begin{align*}
	\tau(H_N)_B = \mathfrak{g} \otimes A \oplus \mathcal{Z}/\mathcal{K}_B \oplus H_B \oplus D	
\end{align*}
is a Lie algebra and is in fact an EALA. We also have an automorphism \\
$B : \tau(S_N) \longrightarrow \tau(S_N)$ such that $B(\mathcal{K}) = \mathcal{K}_B$, which thereby proves that $\tau(H_N) \cong \tau(H_B)$.\\
(3.8) We shall now introduce the notion of twisted HEALA. First let us recall from (3.4) that 
\begin{align*}
\tau(S_N) = \mathfrak{g} \otimes A \oplus \mathcal{Z} \oplus S_N,
\end{align*}
which is an EALA and in particular it carries a non-degenerate symmetric bilinear form. Let $\sigma_1, \ldots, \sigma_N$ be $N$ commuting automorphisms of $\mathfrak{g}$ of orders $n_1, \ldots, n_N$. Let $\xi_i$ be the primitive $n_i$-th root of unity. For $r \in \mathbb{Z}^N$, let
\begin{align*}
\mathfrak{g}_r = \{X \in \mathfrak{g} \ | \ \sigma_i X = \xi_i^{r_i}X \ \forall \ 1 \leqslant i \leqslant N \}.
\end{align*} 
Then we have $\mathfrak{g} = \bigoplus_{0 \leqslant r_i < n_i} \mathfrak{g}_r$. Consider the multiloop algebra
\begin{align*}
L(\mathfrak{g}, \sigma) = \bigoplus_{r \in \mathbb{Z}^N} \mathfrak{g}_{r} \otimes \mathbb{C}t^r.
\end{align*}
We now assume that is a Lie torus (see \cite[Definition 2.2]{ESB}). Consider $A_N = \mathbb{C}[t_1^{\pm n_1}, \ldots, t_N^{\pm n_N}]$ and put $\mathcal{Z}_{\sigma} = \Omega_{A_N}/d_{A_N}, \ \Gamma = \bigoplus_{1 \leqslant i \leqslant N} n_i \mathbb{Z}$. Define
\begin{align*}
\tau(\sigma) = L(\mathfrak{g}, \sigma) \oplus \mathcal{Z}_{\sigma} \oplus S_{\sigma},
\end{align*}
which is called the twisted toroidal EALA relative to $\mathfrak{g}$ and $\sigma = (\sigma_1, \ldots, \sigma_N)$. Here $S_{\sigma} = \text{span} \{D(u,r) \ | \ (u,r) = 0, \ r \in \Gamma \}$. Then clearly $\tau(\sigma) \subseteq \tau(S_N)$ which permits us to consider the restriction of the bilinear form on $\tau(\sigma)$. It is known to be an EALA (see \cite[Section 2]{ESB}).\\
Let us take $H_N(\sigma) = \{D(\overline{s}, s) \ | \ s \in \Gamma \}$ and $\widetilde{H_N}(\sigma) = H_N(\sigma) \oplus D$. Now consider the subalgebra
$L(\mathfrak{g}, \sigma) \oplus \mathcal{Z}_{\sigma} \oplus \widetilde{H_N}(\sigma)$ and further put \\ $\mathcal{K}(\sigma) = \text{span} \{K(u,s) \ | \ u \in \mathbb{C}^N, \ s \in \Gamma, \ (\overline{u},s) = 0 \}$. It can be verified that $[ \widetilde{H_N}(\sigma), \mathcal{K}(\sigma)] \subseteq \mathcal{K}(\sigma)$. Finally define
\begin{align*}
\tau_N (H_N(\sigma)) = L(\mathfrak{g}, \sigma) \oplus \mathcal{Z}_{\sigma}/ \mathcal{K}(\sigma) \oplus \widetilde{H_N}(\sigma) 
\end{align*}
and check that it is indeed an EALA. We refer to this algebra as the twisted Hamiltonian extended affine Lie algebra (twisted HEALA).

\section{The Case $N=2$} \label{S4}
In this section, we assume $N = 2$ and classify all irreducible integrable modules for $\tau(H_2)$ having finite-dimensional weight spaces with respect to $\widetilde{\mathfrak{h}}$ and of non-zero level. The level zero case will be dealt with in Section \ref{S7}.\\
(4.1) Let us first record some special properties for $N = 2$. Recall that $\overline{r} = (r_2, -r_1)$ if $r = (r_1, r_2) \in \mathbb{Z}^2$.\\
(4.1)(1) Suppose $(u,\overline{r}) = 0$. Then it follows that $u$ must be a multiple of $\overline{r}$.\\
(4.1)(2) $\tau(H_2) = \tau(S_2)$.\\
(4.1)(3) $\mathcal{K} = \mathcal{K}_B = 0$ for $B \in GL(2, \mathbb{Z})$.\\
(4.1)(4) $H_B = H_2$ and $((B^T)^{-1}u, Br) = (u,r)$.\\
(4.1)(5) Any $B \in GL(2, \mathbb{Z})$ leaves $\tau(H_2)$ invariant.\\
(4.2) We now define a certain subalgebra of $\tau(H_2)$ and then we shall define a triangular decomposition of $\tau(H_2)$.\\
$\tau(H_2)_{++} = \text{span} \{X(s), \ D(\overline{s},s), \ K(\overline{s},s) \ | \ X \in \mathfrak{g}, \ s=(s_1,s_2) \in \mathbb{Z}^2, \ s_1 > s_2 \} \ \text{and} \\
\tau(H_2)_{+} = \text{span} \{ X_{\alpha}(s), \ D(\overline{s},s), \ K(\overline{s},s) \ | \ X_{\alpha} \in \mathfrak{g}_{\alpha}, \ s=(s_1,s_1) \in \mathbb{Z}^2, \ s_1 > 0 \ \text{or} \ s_1=0,  \ \alpha > 0 \}$.\\
Similarly we define \\
$\tau(H_2)_{--} \ \text{with} \ s_1 < s_2$ and $\tau(H_2)_{-} \ \text{with} \ s_1 = s_2 < 0 \ \text{or} \ s_1=s_2=0, \ \alpha < 0$, \\
alongside $\tau(H_2)_0 = \text{span} \{h, d_1, d_2, K_1, K_2 \ | \ h \in \mathfrak{h} \} = \widetilde{\mathfrak{h}}$. Note that\\
(4.2)(1) $[\tau(H_2)_{++}, \tau(H_2)_{-} \oplus \tau(H_2)_{+}] \subseteq \tau(H_2)_{++}$, \\
(4.2)(2) $[\tau(H_2)_{--}, \tau(H_2)_{-} \oplus \tau(H_2)_{+}] \subseteq \tau(H_2)_{--}$.\\
(4.3) We now define a triangular decomposition for $\tau(H_2)$. Observe that all the above subalgebras are invariant under $\tau(H_2)_0 = \widetilde{\mathfrak{h}}$. Let 
\begin{align*}
\tau(H_2)^+ = \tau(H_2)_{++} \oplus \tau(H_2)_{+} \ , \\
\tau(H_2)^- = \tau(H_2)_{--} \oplus \tau(H_2)_{-} \ . 
\end{align*} 
Then $\tau(H_2) = \tau(H_2)^- \oplus \widetilde{\mathfrak{h}} \oplus \tau(H_2)^+$ is a triangular decomposition of $\tau(H_2)$.\\
(4.4) We now define highest weight irreducible module for $\tau(H_2)$. Take any $\lambda \in \widetilde{\mathfrak{h}}^*$ and let $\mathbb{C}v$ be a $1$-dimensional representation where $\widetilde{\mathfrak{h}}$ acts via $\lambda$. Let $T = \tau(H_2)^{+}$ which acts trivially on $\mathbb{C}v$. Consider the Verma module
\begin{align*}
M(\lambda) = U(\tau(H_2)) \bigotimes_{T} \mathbb{C}v.
\end{align*}      
By standard arguments, it can be easily deduced that $M(\lambda)$ has a unique irreducible quotient $V(\lambda)$.

We now recall the following proposition from \cite[Section 2]{E}. There the module is assumed to be irreducible, but it is not needed in the proof. The proof depends on the root system and the weights being invariant under the Weyl group $W$. The following proposition holds good for any Lie algebra having the same root system as that of the toroidal Lie algebra $\tau$. Hence it also holds good for $\tau(H_2)$.

\bppsn (\cite[Proposition 2.4]{E}, \cite{EJ}) \label{P4.1}
Let $V$ be an integrable module for $\tau(H_2)$ having finite-dimensional weight spaces with respect to $\widetilde{\mathfrak{h}}$. If the central elements $K_1, K_2$ act as a positive integer and zero respectively, then there exists a weight space $V_{\lambda}$ such that\\
(4.5)(1) $X_{\alpha}(r)V_{\lambda} = 0$, \\
(4.5)(2) $D(\overline{r},r)V_{\lambda} = 0$, \\
(4.5)(3) $K(\overline{r},r)V_{\lambda} = 0$, where $r = (r_1,r_2) \in \mathbb{Z}^2$ and $r_1 > 0$.
\eppsn

\noindent (4.6) We now consider the subalgebra $\widetilde{\mathfrak{g}}$ of $\tau(H_2)$ spanned by the elements 
\begin{align*}
\{X(r), \ D(\overline{r},r), \ K_0 = K_1 +K_2, \ d_0 = d_1 + d_2 \ | \ r = (r_1,r_1) \in \mathbb{Z}^2 \}. 
\end{align*}
The Lie bracket is given by 
\begin{align*}
[X(r), Y(s)] = [X,Y](r+s) + r_1 <X,Y> \delta_{r+s,0}(K_1 + K_2). 
\end{align*} 
All the other brackets are obvious. Note that the root system of $\widetilde{\mathfrak{g}}$ is an affine root system given by $\{\alpha + r_1(\delta_1 + \delta_2) \ | \ \alpha \in \Delta \cup \{0 \}, \ r_1 \in \mathbb{Z} \}$. We now recall the following proposition from \cite[Section 2]{E}. 
\bppsn (\cite[Proposition 2.4]{E}) \label{P4.2}
Suppose that $V$ is an integrable module over $\widetilde{\mathfrak{g}}$ having finite-dimensional weight spaces with respect to $\text{span} \{d_1 + d_2, \ h \ | \ h \in \mathfrak{h} \}$. We assume that $K_1 + K_2$ acts by a positive integer. Then there exists a weight space $V_{\lambda}$ of $V$ such that \\
$D(\overline{r},r)V_{\lambda} = 0 \ \forall \ r = (r_1,r_1) \in \mathbb{Z}^2$ and $r_1 > 0$, \\
 $X_{\alpha}(r)V_{\lambda} = 0 \ \forall \ r = (r_1,r_1) \in \mathbb{Z}^2$ with $r_1 > 0$ or $r_1 = 0, \ \alpha > 0$.    
\eppsn

\brmk
Propositions \ref{P4.1} and \ref{P4.2} hold good as they depend on the root system which are equal to the toroidal Lie algebra $\tau$ and $\tau(H_2)$.
\ermk

\bthm \label{T4.4}
Let $V$ be an irreducible integrable module for $\tau(H_2)$ having finite-dimensional weight spaces with respect to $\widetilde{\mathfrak{h}}$. If the level of $V$ is non-zero, then there exists a weight space $V_{\lambda}$ of $V$ such that $\tau(H_2)^{+}.V_{\lambda} = 0$, up to a twist of an automorphism. Hence $V$ is a highest weight module.
\ethm

\begin{proof}
It is a standard fact that if the central operators $K_1$ and $K_2$ act as integers on an integrable modules, then up to a twist of an automorphism, we can assume that $K_1$ acts as $C_1 > 0$ and $K_2$ acts trivially. Now from Proposition \ref{P4.1}, it follows that there exists weight space $V_{\lambda}$ of $V$ such that \\
(4.7) $V_{\lambda + \eta + \delta_r} = 0$ where $\eta \in Q$ and $r_1 > 0$.
Now consider the matrix 
\[B=
\begin{bmatrix}
	a & 1 \\
	a-1 & 1
\end{bmatrix}
\] 
in $GL(2, \mathbb{Z})$. We assume that $2a - 1 > 0$.\\
\textbf{Claim 1.} After twisting the module by the above automorphism $B$, we have $V_{\lambda + \eta + s_1\delta_1 + s_2 \delta_2} = 0 \ \text{where} \ \eta \in Q \ \text{and} \ s_1-s_2 > 0$.\\
Note that
\[B^{-1} = \begin{bmatrix}
	1 & -1 \\
	1-a & a
\end{bmatrix} 
 \ \text{and} \ 
B^{-1}\begin{bmatrix}
	s_1 \\
	s_2
\end{bmatrix} = \begin{bmatrix} 
	s_1 -s_2 \\
	(1-a)s_1 + as_2
\end{bmatrix}. 
\] 
The claim is now obvious from (4.5).\\
Next consider the Lie algebra $\widetilde{\mathfrak{g}}$ where $K_0 = K_1 + K_2$ acts by $(2a-1) > 0$. We can take $d_0 = d_1 + d_2$. Consider the $\widetilde{\mathfrak{g}}$-module $W = \bigoplus_{\eta \in \mathbb{Q}, r \in \mathbb{Z}} V_{\lambda + \eta + r(\delta_1 + \delta_2)}$, which is integrable for $\widetilde{\mathfrak{g}}$ and has finite-dimensional weight spaces. Now by Proposition \ref{P4.2}, there exists a weight space $V_{\mu}$ of $W$ with 
$X_{\alpha} \otimes t_1^{r_1}t_2^{r_1}.V_{\mu} = 0 \ \text{for} \ r_1 > 0 \ \text{or} \ r_1 = 0, \ \alpha > 0$ and $D(\overline{r},r).V_{\mu} = 0 \ \forall \ r = (r_1,r_1) \in \mathbb{Z}^2, \ r_1 > 0$. In order to prove this theorem, it is enough to show the following.\\
\textbf{Claim 2.} $\tau(H_2)_{++}V_{\mu} = 0$.
Let $w \in V_{\mu}$ be arbitrary. Then for any $v \in V_{\lambda}$, there exists $X \in U(\tau(H_2))$ such that $Xv=w$. By the PBW theorem, there exist $X_1 \in U(\tau(H_2)_{--}),$ $X_2 \in U(\tau(H_2)_{-}),$ $X_3 \in U(\tau(H_2)_{+})$ and $X_{4} \in U(\tau(H_2)_{++})$ such that $X = X_1X_2X_3X_4$. By Claim 1, $X_{\alpha}v = 0$ unless $X_{\alpha}$ is a scalar. Note that $\mu - \lambda = \eta +  r_1 (\delta_1 +\delta_2)$ for some $r_1 \in \mathbb{Z}$ and $\eta \in Q$. Now by weight arguments, $X_1$ has to be a scalar. Thus we have $X_2X_3v = w$.\\
\textbf{Claim 3.} $\tau(H_2)_{++}.X_2X_3v = 0$.\\
The claim follows from (4.2)(1). This concludes the proof of Claim 2.
\end{proof} 

\brmk
These non-zero level integrable modules have been classified in \cite{CLT1}. But their triangular decomposition is different. For them, the middle part contains some special subalgebra of $DerA$ and needs special attention. In our case, the middle part is just $\widetilde{\mathfrak{h}}$.
\ermk 

\brmk
\
\begin{enumerate}
\item Compare the modules in \cite{CLT1} with the highest weight modules defined in this section.
\item Classify $\lambda \in \widetilde{\mathfrak{h}}^*$ for which $V(\lambda)$ has finite-dimensional weight spaces.
\item Classify $\lambda \in \widetilde{\mathfrak{h}}^*$ for which $V(\lambda)$ is integrable.   
\end{enumerate}
\ermk

\section{Classification of irreducible integrable modules for $\tau(H_N)$} \label{S5}
In this section and in the next section, we classify irreducible integrable modules for $\tau(H_N)$ having finite-dimensional weight spaces with respect to $\widetilde{\mathfrak{h}}$. We assume that $N = 2m \geqslant 4$.\\ 
(5.1) We first define a triangular decomposition for $\tau(H_N)$ similar to the earlier section but the zeroth component will be infinite-dimensional.\\
$\tau(H_N)_{++} = \text{span} \{X(r), \ D(\overline{r},r), \ K(\overline{r},r) \ | \ X \in \mathfrak{g}, \ r \in \mathbb{Z}^N, \ r_m > r_{2m} \}, \\
\tau(H_N)_{+} = \text{span} \{ X_{\alpha}(r), \ D(\overline{r},r), \ K(\overline{r},r) \ | \ X_{\alpha} \in \mathfrak{g}_{\alpha}, \ r \in \mathbb{Z}^N, \ r_m = r_{2m} > 0 \ \text{or} \ r_m = r_{2m}=0, \ \alpha > 0 \}$ and \\
$\tau(H_N)_{0} = \text{span} \{ h(r), \ D(\overline{r},r), \ K(\overline{r},r), \ D(u,0), \ K(u,0) \ | \ h \in \mathfrak{h}, \ r \in \mathbb{Z}^N, \\ r_m = r_{2m} = 0, \ u \in \mathbb{C}^N \}$.\\
Note that $\widetilde{\mathfrak{h}} \subseteq \tau(H_N)_0$. Similarly we define $\tau(H_N)_{--}$ and $\tau(H_N)_{-}$. Observe that all the five algebras are closed under Lie bracket. Let \\
$\tau(H_N)^+ = \tau(H_N)_{++} \oplus \tau(H_N)_{+}$ and $\tau(H_N)^- = \tau(H_N)_{--} \oplus \tau(H_N)_{-}$.\\
(5.2) Then $\tau(H_N) = \tau(H_N)^- \oplus \tau(H_N)_{0} \oplus \tau(H_N)^+$ is clearly a triangular decomposition of $\tau(H_N)$. Note that $\tau(H_N)_{++},  \tau(H_N)_{+},  \tau(H_N)_{-}, \tau(H_N)_{--}$ are invariant under $\tau(H_N)_{0}$ and $[\tau(H_N)_{++}, \tau(H_N)_{-} \oplus \tau(H_N)_{+}] \subseteq \tau(H_N)_{++}$.\\
(5.3) We now define highest weight modules for $\tau(H_N)$. Let $U^{+} = U(\tau(H_N)^{+})$ and $U^{-} = U(\tau(H_N)^{-})$. Suppose that $W$ is an irreducible module for $\tau(H_N)_{0}$ having finite-dimensional weight spaces with respect to $\widetilde{\mathfrak{h}}$. Consider the Verma module given by 
\begin{align*}
M(W) = U \bigotimes_{U^+ \oplus U^0} W, 
\end{align*}           
where $U = U(\tau(H_N))$ and $U^+$ acts trivially on $W$. By standard arguments, $M(W)$ has a unique irreducible quotient which we shall denote by $V(W)$.\\
(5.4) $r \in \mathbb{Z}^{N-2}$ is treated as a vector in $\mathbb{Z}^{N}$ by setting $r_m = r_{2m} = 0$. We shall follow a similar notation for a vector in $\mathbb{C}^{N-2}$.\\
(5.5) We now fix an irreducible integrable module $V$ for $\tau(H_N)$ having finite-dimensional weight spaces with respect to $\widetilde{\mathfrak{h}}$. We assume that each $K_i$ acts trivially except for $i=m$. We assume that $K_m$ acts by a positive integer.

\brmk
In the non-zero level case, we can always assume (5.5) up to a twist of an automorphism B (see \cite{E}). But then we need to work with $\tau(H_B)$. To avoid notational complexity, we assume (5.5).
\ermk 

\bppsn
Let $V$ be as in (5.5). Then there exists a weight space $V_{\lambda}$ of $V$ such that $V_{\lambda + \eta + \delta_r} = 0$ for $r_m > 0, \ \eta \in Q$.
\eppsn 

\begin{proof}
This follows from Proposition \ref{P4.1} wich is valid for several variables (see \cite{E}).
\end{proof}

We shall now twist the module by an automorphism $B \in GL(N, \mathbb{Z})$. Define $B = (b_{ij})_{N \times N}$, where $b_{ii} = 1$ for  $i \neq m, 2m$ and $b_{ij} = 0$ for $i \neq j$ except $i =m$ and $j = 2m$ or $i=2m$ and $j = m$. That leaves the following matrix to be defined.
\[
\begin{bmatrix}
	b_{m,m} & b_{m,2m} \\
	b_{2m,m} & b_{2m,2m}
\end{bmatrix}
\] 	
which we assume to be
\[
\begin{bmatrix}
a & 1 \\
a-1 & 1
\end{bmatrix} \ \text{with} \ 2a-1 > 0
\]
Notice that the inverse is  
\[ \begin{bmatrix}
	1 & -1 \\
	1-a & a
\end{bmatrix} 
\ \text{and} \ 
\begin{bmatrix}
	1 & -1 \\
	1-a & a
\end{bmatrix}\begin{bmatrix}
	s_m \\
	s_{2m}
\end{bmatrix} = \begin{bmatrix} 
	s_m -s_{2m} \\
	(1-a)s_m + as_{2m}
\end{bmatrix}. 
\] 
Now it is obvious that after twisting the module by $B$, we have a weight space $V_{\lambda}$ of $V$ such that\\
(5.6) $ V_{\lambda + \eta + \delta_s} = 0, \ \text{where} \ \eta \in Q, \ s_m - s_{2m} > 0$.
\bthm
Let $V$ be as above. Then up to an automorphism, there exists a weight space $V_{\mu}$ of $V$ such that $\tau(H_N)^+.V_{\mu} = 0$.
\ethm

\begin{proof}
Let $\widetilde{\mathfrak{g}} = \text{span} \{X(r), \ D(\overline{r},r), \ K(\overline{r},r), \ D(u,0), \ K(u,0) \ | \ X \in \mathfrak{g}, \ r \in \mathbb{Z}^N, \ r_m = r_{2m}, \ u \in \mathbb{C}^N \}$, which is a Lie subalgebra of $\tau(H_N)$. Consider the $\widetilde{\mathfrak{g}}$-module
\begin{align*}
\bigoplus_{\eta \in Q,\ l \in \mathbb{Z}^N, \ l_m=l_{2m}} V_{\lambda + \eta + \delta_l}.
\end{align*} 
Now by Proposition \ref{P4.2} (see \cite[Proposition 2.4]{E} for a general version), there exists $\mu = \lambda + \eta^{\prime} + \delta_r$ where $\eta^{\prime} \in Q, \ r_m = r_{2m}$ such that $V_{\mu + \eta + \delta_s} = 0$ for $\eta \in Q, \ s \in \mathbb{Z}^N, \ s_m=s_{2m} > 0$ or $s_m = s_{2m} = 0$ but $\eta \in Q^+$.\\
\textbf{Claim.} $\tau(H_N)_{++}V_{\mu} = 0$.\\
Let $w \in V_{\mu}$ and $v \in V_{\lambda}$. Then by the irreducibility of $V$, there exists $X \in U(\tau(H_N))$ such that $Xv=w$. But by PBW theorem, there exist $X_1 \in U(\tau(H_N)_{--}), \ X_2 \in U(\tau(H_N)_{-}), \ X_3 \in U(\tau(H_N)_{+}), \ X_4 \in U(\tau(H_N)_{++})$ and $X_0 \in U(\tau(H_N)_{0})$ such that $X = X_1X_2X_3X_4X_0$. But by (5.6), $X_4X_0v = 0$ unless $X_4$ is a scalar. Note that $\mu - \lambda = \eta + \delta_l$ for $\eta \in Q$ and $l_m = l_{2m}$. Now by weight reasons, $X_1$ has to be a scalar. Thus $w = X_2X_3X_0v$. Subsequently using (5.2), it is easy to see that $\tau(H_N)_{++}w = 0$. This completes the proof of this theorem.  
\end{proof}

Now consider $W = \{v \in V \ | \ \tau(H_N)^+v = 0 \}$. From above, $W$ is non-zero. In fact, $W$ is an irreducible module for $\tau(H_N)_0$. Thus $V \cong V(W)$ as defined in (5.3). We shall describe $W$ in more detail in the next section using results from Appendix and \cite{T}.

\brmk
One can take another natural triangular decomposition by declaring the positive part with $r_m > 0$. Then the zeroth component (defined by $r_m = 0$) will contain a loop algebra of $\mathfrak{h} \otimes A_{n-2} \oplus H_{N-2}$. We do not have classification result for such loop algebras. We except irreducible cuspidal modules for such loop modules to be evaluation modules at a single point. We now explain the loop. Consider $\overline{H_N} = \text{span} \{h_r \in H_N \ | \ r_m = 0 \}$. For $r \in \mathbb{Z}^N$ with $r_m = 0$, let $r^{\prime} = (r_1,r_2, \ldots, r_{m-1}, 0, r_{m+1}, \ldots, 0)$. Henceforth define $D(\overline{r}^{\prime}, r^{\prime}) \otimes t^{r_{2m}} = D(\overline{r}, r) \ (r_m = 0)$. It can be verified that for $r_m = 0 =s_m$, we have $[D(\overline{r}^{\prime}, r^{\prime}) \otimes t^{r_{2m}}, D(\overline{s}^{\prime}, s^{\prime}) \otimes t^{s_{2m}}] = [D(\overline{r}, r), D(\overline{s}, s)] = [D(\overline{r}^{\prime}, r^{\prime}), D(\overline{s}^{\prime}, s^{\prime})] \otimes t^{r_{2m} + s_{2m}} = D(\overline{r+s}^{\prime}, (r+s)^{\prime}) \otimes  t^{r_{2m} + s_{2m}}$. Note that for $r,s \in \mathbb{Z}^N$ with $r_m = s_m =0$, $(\overline{r}, s) = (\overline{r}^{\prime}, s^{\prime})$. So as we do not have classification result for the loop algebra of $\mathfrak{h} \otimes A_{N-2} \oplus H_N$, we considered another triangular decomposition which is defined at the start of this section. In this case, we do have a classification result for $\mathfrak{h} \otimes A_{N-2} \oplus H_{N-2}$ due to \cite{T}, which is stated in Theorem \ref{T6.4}.    
\ermk

\section{Description of highest weight space} \label{S6}
In this section, we describe the $\tau(H_N)_{0}$-module $W$ in more detail by using results from \cite{T} and our Appendix. It is a weight module for $\widetilde{\mathfrak{h}}$ where $K_i, i \notin \{m, 2m \}$ acts trivially. $d_m$ and $d_{2m}$ act as scalars, say $g_m$ and $g_{2m}$. Note that $\mathfrak{h}$ is central in $\tau(H_N)_0$ and hence acts by a single linear functional, say $\overline{\lambda}$. Since the module is integrable, we see that $\overline{\lambda}$ is dominant integral. We shall describe $W$ when $\overline{\lambda}$ is non-zero.\\
(6.1) Let $\lambda = \overline{\lambda} + g_m \delta_m + g_{2m} \delta_{2m} + aKw_m + (a-1)Kw_{2m}$. Note that $K_m$ acts as $aK$ and $K_{2m}$ acts as $(a-1)K$ ($K$ is as defined in Section \ref{S5} which is the $m$-th co-ordinate of central element).\\
Define $\lambda_r = \lambda + \delta_r, \ r \in \mathbb{Z}^{N-2} (r_m = 0 = r_{2m})$. Let $W = \oplus_{\lambda \in \widetilde{\mathfrak{h}}^*} W_{\lambda}$ and $P(W) = \{\lambda \in \widetilde{\mathfrak{h}}^* \ | \ W_{\lambda} \neq 0 \}$. Then it follows that $P(W) \subseteq \{\lambda_r \ | \ r \in \mathbb{Z}^{N-2} \}$.\\
(6.2) We shall first prove that the weight spaces are uniformly bounded for $\overline{\lambda} \neq 0$. We also show that $\tau(H_N)_0 \cap \mathcal{Z}/\mathcal{K}$ acts trivially on $W$.

\bppsn \label{P6.1}
Suppose that $\overline{\lambda} \neq 0$. Then there exists an integer $p$ such that dim $W_{\lambda} \leqslant p \ \forall \ \lambda \in P(W)$.
\eppsn

\begin{proof}
We shall use the fact that $K_i$ acts trivially on $W$ for $i \notin \{m, 2m \}$ and that $P(W)$ is invariant under the Weyl group. By using a similar argument as in \cite[Lemma 3.1]{CG} (see also \cite{ESB}), it follows that there exists $s_{\overline{\lambda}} > 0$ depending on $\overline{\lambda}$ (hence unique) with the following property : For any $\lambda_r \in P(W)$, there exists $w_i \in W$ such that $w_i(\lambda_r) = \lambda_s$ where $0 \leqslant s_i < s_{\overline{\lambda}}$ for $i \notin \{m,2m \}$.\\
Put $p =$ max$\{ \text{dim} W_{\lambda_s} \ | \ 0 \leqslant s_i < s_{\overline{\lambda}} \}$. Thus dim $W_{\lambda_r} \leqslant p \ \forall \ \lambda_r \in P(W)$.
\end{proof}   

\bppsn \label{P6.2}
$\mathcal{Z}^{\prime} = U(H_N)_0 \cap \mathcal{Z}/\mathcal{K}$ acts trivially on $W$.
\eppsn 

\begin{proof}
First consider $u_mt^rK_m + u_{2m}t^rK_{2m}, \ 0 \neq r \in \mathbb{Z}^{N-2}, \ u_m, u_{2m} \in \mathbb{C}$. From the definition of $\mathcal{K}$, it follows that they are all zero vectors. Thus $\mathcal{Z}^{\prime}$ consists of $K(u,r), u \in \mathbb{C}^{N-2}$ and $r \in \mathbb{Z}^{N-2}$. Let $M = \oplus_{0 \leqslant s_i < s_{\overline{\lambda}}} W_{\lambda_s}$, where we have dim $M < \infty$. Let $\tau_{N-2} = \mathfrak{g} \otimes A_{N-2} \oplus \mathcal{Z}/\mathcal{K} \oplus \text{span} \{d_i \ | \ i \notin \{m, 2m \} \}$. Now note that $\tau_{N-2}$ is a quotient of $(N-2)$ variables toroidal Lie algebra as we are quotienting out some part of the center. Consider $W^{\prime} = U(\tau_{N-2})W$ which is a $\tau_{N-2}$-module. Suppose that there exists a decreasing sequence of $\tau_{N-2}$-submodules given by
\begin{align*}
\ldots \subseteq W_{i+1} \subseteq W_i \subseteq \ldots \subseteq W^{\prime},
\end{align*}
such that $W_i \cap W \neq 0$. Then it follows that $M_i = W_i \cap W \neq 0$. Thus we have a decreasing sequence of subspaces of $M$. Since $M$ is finite-dimensional, the above chain of submodules must necessarily terminate. Thus there exists a $\tau_{N-2}$-submodule $W^{\prime \prime}$ of $W^{\prime}$ such that $W^{\prime \prime} \cap W$ is an irreducible module for $\mathfrak{h} \otimes A_{N-2} \oplus \mathcal{Z}^{\prime} \oplus \text{span} \{d_i \ | \ i \notin \{m, 2m \} \}$. This $W^{\prime \prime}$ need not be irreducible, but one can take the unique irreducible quotient, where $W^{\prime \prime} \cap W$ goes injectively to the quotient. Thus we have an irreducible module for $\tau_{N-2}$ where the zero degree center acts trivially. Then by \cite[Proposition 4.13]{E}, it follows that $\mathcal{Z}^{\prime}$ also acts trivially. Consider $W^{\prime \prime \prime} = \{w \in W \ | \ \mathcal{Z}^{\prime}w = 0 \}$, which is a non-zero $\tau(H_N)_0$-submodule of $W$. But since $W$ is irreducible, it follows that $W^{\prime \prime \prime} = W$. Thus the proposition is proved.  
\end{proof}

\bppsn \label{P6.3}
Let $\overline{\lambda} = 0$. Then $\mathfrak{h} \otimes A_{N-1} \oplus \mathcal{Z}^{\prime}$ acts trivially on $W$.
\eppsn

\begin{proof}
Let $X_{\alpha} \in \mathfrak{g}_{\alpha}$ and $Y_{\alpha} \in \mathfrak{g}_{-\alpha}$ with $[X_{\alpha}, Y_{\alpha}] = h_{\alpha}, \ [h_{\alpha}, X_{\alpha}] = 2X_{\alpha}$, $[h_{\alpha}, Y_{\alpha}] = -2Y_{\alpha}$ so that $\text{span} \{X_{\alpha}, Y_{\alpha}, h_{\alpha} \} \cong \mathfrak{sl}_2$. \\
\textbf{Claim.} $h_{\alpha}(r)$ acts trivially on $W$ for all  $r \neq 0$ if and only if $\overline{\lambda}(h_{\alpha}) = 0$.\\
Consider $\text{span} \{X_{\alpha}(r), Y_{\alpha}(-r), h_{\alpha} + \sum_{i \notin \{m,2m \}} r_iK_i \} \cong \mathfrak{sl}_2$. Suppose that $\overline{\lambda}(h_{\alpha}) = 0$. Then using integrability, it follows that $Y_{\alpha}(-r)$ acts trivially on $W$. Thereby considering $[X_{\alpha}(r), Y_{-\alpha}(s)] = h_{\alpha} \otimes t^{r+s} + K(r,r+s)$, we see that $h_{\alpha} \otimes t^{r+s}$ acts trivially.\\
Conversely suppose that $h_{\alpha}(r)$ acts trivially for all $r \neq 0$. Then consider 
$h_{\alpha}(r+s) = X_{\alpha}(r)Y_{\alpha}(s) - Y_{\alpha}(s)X_{\alpha}(r)$, which immediately implies that $X_{\alpha}(r)Y_{\alpha}(s)w = 0$. By integrable theory of $\mathfrak{sl}_2$, it follows that $Y_{\alpha}(s)$ is a highest weight vector and hence $\overline{\lambda}(h_{\alpha}) = 0$.    
\end{proof}

\noindent (6.3) In this subsection, we record a very important classification theorem for Jet modules over $\widetilde{H_N}$ from \cite{T}.

A module $W_N$ for $\widetilde{H_N}$ is said to be a Jet module if $A_N = A$ acts on $W_N$ associatively, in the sense $t^rt^s = t^{r+s}$ on $W_N$ and $t^0$ acts as the identity operator. In a Jet module, $t^r$ acts injectively for all $r \in \mathbb{Z}^N$. Hence all the weight spaces have the same dimension. So we can write $W_N = V_N \otimes A_N$ for some vector space $V_N$. It is also clear that the weights of $W_N$ are of the form $\lambda + \mathbb{Z}^N$ for some $\lambda \in \mathbb{C}^N$ ($N = 2m$), provided $W_N$ is indecomposable.

Let $E_{ij}$ be the $N \times N$ matrix with $(i,j)$-th entry one and zeros elsewhere. Then the symplectic Lie algebra $\mathfrak{sp}_N$ is spanned by the $2N^2 - N$ vectors
\begin{align*}
E_{ij} - E_{m+j,m+i}, \ E_{i,m+j} + E_{j,m+i} \ \text{and} \ E_{i,m+i}, \  E_{m+i,i} \ \text{for} \ 1 \leqslant i,j \leqslant N.
\end{align*}

\bthm \cite[Theorem 5.2]{T} \label{T6.4} 
Let $W_N$ be an irreducible Jet module for $\widetilde{H_N}$ having finite-dimensional weight spaces with respect to $\widetilde{\mathfrak{h}}$. Then there exist a finite-dimensional irreducible module $V_N$ for $\mathfrak{sp}_N$ and $u, w \in \mathbb{C}^N$ such that $W_N = V_N \otimes A_N$ with the following actions.\\
(1) $d_i (v \otimes t^s) = (s_i + u_i)(v \otimes t^s)$, \\
(2) $h(r) (v \otimes t^s) = (\overline{r}, s) (v \otimes t^{r+s}) + ( \sum_{i} r_{m+i}w_{i+m} - \sum_{i}r_iw_i)(v \otimes t^{r+s}) + \sum_{i} \big( (r_{m+i}^2 E_{m+i,i}) + r_im_{i+m}(E_{ii} - E_{m+i,m+i}) + (r_i^2E_{i,m+i}) \big ) (v \otimes t^{r+s}) + \\ \sum_{i,j} \big (r_{m+i}r_{m+j}(E_{m+j,i} + E_{m+i,j}) + r_ir_{m+j} (E_{ij} - E_{m+i,m+j}) - \\ r_ir_j(E_{i,m+j} + E_{j,m+i}) \big ) (v \otimes t^{r+s})$, \\
(3) $t^r (v \otimes t^s) = v \otimes t^{r+s}$.
\ethm  

\noindent (6.4) Now we are in a position to describe the $\tau(H_N)_0$-module $W$ with the help of Theorem \ref{T6.4} and Appendix. We already know that $W$ ia a module over $\mathfrak{h} \otimes A_{N-2} \oplus \widetilde{H_{N-2}}$ (see Proposition \ref{P6.2}). We assume that $\mathfrak{h} \otimes A_{N-2}$ acts non-trivially. Let $\alpha \in \Delta$ and suppose that $\lambda(h_{\alpha}) \neq 0$. Then we know that $h_{\alpha} \otimes t^r$ ($r \neq 0$) acts non-trivially (by Claim in Proposition \ref{P6.3}). Again by Theorem \ref{T9.1} (in Appendix), there exist non-zero scalars $\lambda_{\alpha}, \mu_{\alpha}$ such that $\lambda_{\alpha}^2 = \mu_{\alpha} \lambda(h_{\alpha})$ satisfy $(h_{\alpha} \otimes t^r)(h_{\alpha} \otimes t^s) = \lambda_{\alpha} (h_{\alpha} \otimes t^{r+s}), \ r,s \neq 0, r+s \neq 0$.\\
$(h_{\alpha} \otimes t^r)(h_{\alpha} \otimes t^{-r}) = \mu_{\alpha}h_{\alpha} = \mu_{\alpha} \lambda(h_{\alpha})$. We can now take 
\begin{align*}
A = \bigg \{\dfrac{h_{\alpha}}{\lambda(h_{\alpha})}, \ \dfrac{h_{\alpha} \otimes t^r}{\lambda_{\alpha}} \ \big | \ r \neq 0  \bigg \}
\end{align*}
and the action of $A$ is associative. Hence $W$ becomes a Jet module for $\widetilde{H_N}$ and thus $W \cong W_N$ (by Theorem \ref{T6.4}). Consequently we have an irreducible module $V(W_N)$ over $\tau(H_N)$ (see (5.3)). Let us now state the main theorem of this section which follows from earlier results.

\bthm
Let $V$ be an irreducible integrable module of non-zero level for $\tau(H_N)$ having finite-dimensional weight spaces with respect to $\widetilde{\mathfrak{h}}$. Let $W$ be the corresponding highest weight space (up to an automorphism). Assume that $\mathfrak{h}$ acts non-trivially on $W$. Then $V \cong V(W_{N-2})$ as $\tau(H_N)$-modules, up to a twist of an automorphism. 
\ethm

\section{The level zero case}\label{S7}
In this section, we classify the level zero irreducible integrable modules for $\tau(H_N)$ with finite-dimensional weight spaces. The strategy is to prove that such a module is a highest weight module with respect to a suitable triangular decomposition. We prove that $\mathcal{Z}/\mathcal{K}$ acts trivially on the highest weight space (in fact, it acts trivially on the whole space). Next we prove that $\mathfrak{h} \otimes A$ acts ``associatively" using results from the Appendix. Then using the classification result from \cite{T}, we completely describe the highest weight space, thereby showing that our module is the irreducible quotient of the induced module. Throughout this section, we always assume that $\mathfrak{g} \otimes A$ acts non-trivially on our irreducible module. At the end, we provide a realization of such modules. \\
(7.1) Throughout this section, we shall work with the following triangular decomposition of $\tau(H_N)$.\\
$\tau(H_N)(+) = \text{span} \{X_{\alpha}(r) \ | \ \alpha \in \Delta^+, \ r \in \mathbb{Z}^N \}$, \\
$\tau(H_N)(-) = \text{span} \{X_{\alpha}(r) \ | \ \alpha \in \Delta^-, \ r \in \mathbb{Z}^N \}$,\\
$\tau(H_N)(0) = \text{span} \{h(r), \ K(\overline{r},r), \ K(u,0), \ D(u,0), \ D(\overline{r},r) \ | \ u \in \mathbb{C}^N, \ \ r \in \mathbb{Z}^N, \ h \in \mathfrak{h} \}$.\\
(7.2) We first recall a well-known result from \cite{C} (see \cite{E} for the multivariable case). Let $L(\mathfrak{g}) = \mathfrak{g} \otimes \mathbb{C}[t,t^{-1}] \oplus \mathbb{C}d$ be a centerless affine Lie algebra.

\bthm \cite[Section 4]{C} \label{T7.1}
Let $V$ be an irreducible integrable module for $L(\mathfrak{g})$ having finite-dimensional weight spaces with respect to $\mathfrak{h} \oplus \mathbb{C}d$. Then there exists some $p \in \mathbb{N}$ and non-zero distinct complex numbers $a_1, a_2, \ldots, a_p$ and dominant integral weights $\lambda_1, \ldots, \lambda_p \in \mathfrak{h}^*$ such that $V$ is isomorphic to an irreducible component of $\big(\otimes_{i=1}^{p} V(\lambda_i) \big) \otimes \mathbb{C}[t,t^{-1}]$, where $V(\lambda_i)$ is the finite-dimensional irreducible highest weight module for $\mathfrak{g}$. The $L(\mathfrak{g})$ action is given by the following. 
\begin{align*}
(X \otimes t^r)(v_1 \otimes \ldots \otimes v_p \otimes t^q) = \sum_{i} a_i^rv_1 \otimes \ldots \otimes X.v_i \otimes \ldots \otimes v_p \otimes t^{r+q} \\
\text{for all} \ X \in \mathfrak{g} \ \text{and} \ \ q,r \in \mathbb{Z}.\\
\text{In particular}, \ (h \otimes t^r)(v \otimes t^q) = \sum_{i} \lambda_i(h) a_i^r v \otimes t^{r+q} \ \forall \ h \in \mathfrak{h}, \\
v = v_1 \otimes \ldots \otimes v_p, \ \text{where} \ v_i \ \text{is a highest weight vector for} \ V(\lambda_i).\\
\text{In case} \ \sum_{i} \lambda_i(h) a_i^r \neq 0 \ \forall \ r, \ \text{then} \ \big (\otimes_{i=1}^{p} V(\lambda_i) \big ) \otimes \mathbb{C}[t,t^{-1}] \ \text{is irreducible}.  
\end{align*}
\ethm  

\noindent (7.3) We also recall a result which is essentially due to Chari \cite{C}. It is proved for one variable in \cite{C} but here we need the corresponding result for several variables. We briefly sketch the proof maintaining the same notations that we have utilized in Section \ref{S2}.

\bthm \cite{C,E} \label{T7.2}
Suppose that $V$ is an irreducible integrable module for $\mathfrak{g} \otimes A \oplus \mathcal{Z} \oplus D$ having finite-dimensional weight spaces with respect to $\widetilde{\mathfrak{h}}$ and let each $K_i$ act trivially. Then there exists a weight space $V_{\lambda}$ of $V$ such that
\begin{align*}
	V_{\lambda + \alpha + \delta_r} = 0 \ \forall \ \alpha > 0, \ r \in \mathbb{Z}^N. 
\end{align*}
\ethm

\begin{proof}
From the proof of \cite[Lemma 2.6]{E}, we have a weight space $V_{\mu}$ such that $V_{\mu + \alpha} = 0 \ \forall \ \alpha > 0$. Let $V_{\mu + \alpha + \delta_r} \neq 0$ for some $\alpha > 0$ and some $\delta_r$.\\
\textbf{Claim.} $V_{\lambda + \beta + \delta_s} = 0$ where $\lambda = \mu + \alpha + \delta_r \ \forall \ \beta > 0$ and $s \in \mathbb{Z}^N$.\\
Suppose not for some $\beta > 0$ and $\delta_s$. Since $\alpha, \beta \in \Delta^+$, either $(\alpha + \beta, \alpha) > 0$ or $(\alpha + \beta, \beta) > 0$. Assume that $(\alpha + \beta, \alpha) > 0$ and set $\gamma = \alpha+ \delta_r + \delta_s$. Then $(\alpha + \beta, \gamma) > 0$. Hence by (2.13)(4), we have $V_{\lambda + \beta + \delta_s - \gamma} \neq 0 \implies V_{\mu + s} \neq 0$, a contradiction. Hence the claim. This completes the proof of this theorem. 

\noindent (7.4) For the rest of this section, let us fix an irreducible integrable module for $\tau(H_N)$ having finite-dimensional weight spaces with respect to $\widetilde{\mathfrak{h}}$. Assume that each $K_i$ acts trivially and $\mathfrak{g} \otimes A$ acts non-trivially. Set
\begin{align*}
W^{\circ} = \{v \in V \ | \ \tau(H_N)(+)v = 0 \}.
\end{align*} 
It follows from Theorem \ref{T7.2} that $W^{\circ}$ is non-zero. By PBW theorem, it further follows that $W^{\circ}$ is in fact an irreducible module over $\tau(H_N)(0)$. 

\bppsn
The weight spaces of $W^{\circ}$ are uniformly bounded.
\eppsn
\begin{proof}
This follows from arguments similar to Proposition \ref{P6.1} as each $K_i$ acts trivially.
\end{proof}

\bppsn
$\mathcal{Z}/\mathcal{K}$ acts trivially on $W^{\circ}$.
\eppsn
\begin{proof}
Follows from similar arguments as in Proposition \ref{P6.2}.
\end{proof}

Subsequently $W^{\circ}$ is a module for $\mathfrak{h} \otimes A \oplus \widetilde{H_N}$. It follows that $\mathfrak{h}$ is central and hence acts by a single linear functional, say $\overline{\lambda}$. Obviously $\overline{\lambda}$ is dominant integral. Since we are assuming that $\mathfrak{g} \otimes A$ acts trivially, we see that $\overline{\lambda} \neq 0$ (see the claim in Proposition \ref{P6.3}).

Fix a vector $v$ in $W^{\circ}$ and let $d_iv = g_iv$. Then any weight of $W^{\circ}$ looks like $\overline{\lambda} + \delta_r + \sum_{i} g_i \delta_i$, where $g_1, g_2, \ldots, g_N$ are fixed constants. Let $\alpha \in \Delta^+$ such that $\lambda(h_{\alpha}) \neq 0$. Then by the claim in Proposition \ref{P6.3}, we see that $h_{\alpha} \otimes t^r$ acts non-trivially for each $r \neq 0$. Now by Theorem \ref{T9.1} (appendix), we see that there exist non-zero scalars $\lambda_{\alpha}, \mu_{\alpha}$ such that on $W^{\circ}$, we have\\
(7.4) $(h_{\alpha} \otimes t^r)(h_{\alpha} \otimes t^s) = \lambda_{\alpha}(h_{\alpha} \otimes t^{r+s}), \ r,s \neq 0, \ r+s \neq 0$, \\ 
$(h_{\alpha} \otimes t^r)(h_{\alpha} \otimes t^{-r}) = \mu_{\alpha} \lambda(h_{\alpha})$ and $\lambda_{\alpha}^2 = \mu_{\alpha} \lambda(h_{\alpha})$.

\blmma
$\lambda_{\alpha} = \lambda(h_{\alpha}) \ \forall \ \alpha \in \Delta^{+}$ whenever $\lambda(h_{\alpha}) \neq 0$.
\elmma

\begin{proof}
By Theorem \ref{T6.4}, we have $W^{\circ} = W_N = V_N \otimes A$ as $\widetilde{H_N}$-module. Take $v \in V_N$ and look at $\mathbb{C}v \otimes A_N$, which is invariant under $\mathbb{C}h_{\alpha} \otimes A$ and in fact it is irreducible. Take the centerless affine Lie algebra $\mathfrak{g}_{aff}^1 = \mathfrak{g} \otimes \mathbb{C}[t, t^{-1}] \oplus \mathbb{C}d$ and consider the module 
$U(\mathfrak{g}_{aff}^1)(\mathbb{C}h_{\alpha} \otimes A)(\mathbb{C}v \otimes A) = W^{\prime}$. Define
\begin{align*}
\mathfrak{g}_{aff}^{1,+} = \{\mathfrak{g}_{\alpha} \otimes \mathbb{C}[t,t^{-1}], \ \alpha \in \Delta^{+} \}, \,\,\,\ \mathfrak{g}_{aff}^{1,-} = \{\mathfrak{g}_{\alpha} \otimes \mathbb{C}[t,t^{-1}], \ \alpha \in \Delta^{-} \}.
\end{align*} 
Then clearly $W^{\prime}$ is a highest weight module for $\mathfrak{g}_{aff}^1$. The module $W^{\prime}$ may not be irreducible for $\mathfrak{g}_{aff}^1$, but it has an unique irreducible quotient where $\mathbb{C}v \otimes A$ goes injectively. Consequently by Theorem \ref{T7.1}, there exists $p \in \mathbb{N}$, non-zero distinct complex numbers $a_1, \ldots, a_p$ and dominant integral weights $\lambda_1, \ldots, \lambda_p$ such that $(h_{\alpha} \otimes t^q)(v \otimes t^r) = \sum_{i} \lambda_i(h_{\alpha})a_i^q(v \otimes t^{q + r})$. But by Theorem \ref{T6.4} and (7.4), we have $(h_{\alpha} \otimes t^q)(v \otimes t^r) = \lambda_{\alpha}(v \otimes t^{q + r})$. This implies that $\lambda_{\alpha} = \sum_{i} \lambda_i(h_{\alpha})a_i^s \ \forall \ s \neq 0$. Note that $\lambda_{\alpha}$ is independent of $s$. Consider the Vandermonde matrix given by 
\[M =
\begin{bmatrix}
	a_1 & a_2 & \ldots & a_p \\
	a_1^2 & a_2^2 & \ldots & a_p^2 \\
	\ldots & \ldots & \ldots & \ldots \\
	a_1^p & a_2^p & \ldots & a_p^p
\end{bmatrix}.
\] 	   
\end{proof} 
Let
\[\underline{\lambda} = \begin{bmatrix}
	\lambda_1(h_{\alpha}) \\
	\lambda_2(h_{\alpha}) \\
	\ldots \\
	\lambda_p(h_{\alpha})
\end{bmatrix} 
\ \,\,\,\ \text{and} \,\,\,\ \overline{D} = \ 
\begin{bmatrix}
	a_1 & 0 & \ldots & 0 \\
	0 & a_2 & \ldots & 0 \\
	\ldots & \ldots & \ldots & \ldots \\
	0  & 0 & \ldots & a_p
\end{bmatrix}.
\] 
Then it follows from above that $M \underline{\lambda} = M \overline{D} \underline{\lambda}$. Since $M$ is invertible, we see that $\underline{\lambda} = D \overline{\lambda} \implies \lambda_i(h_{\alpha}) = a_i \lambda_i(h_{\alpha})$. Hence $a_i =1$ and also $p=1$. This proves that $\lambda_{\alpha} = \lambda(h_{\alpha})$. Note that $\lambda_{\alpha} = \mu_{\alpha} = c = \lambda(h_{\alpha})$ as $h_{\alpha} \otimes t^0$ acts as $\lambda(h_{\alpha})$. 
\end{proof} 

\noindent (7.5) We have now described $W^{\circ}$ very explicitly. It is a module for $\mathfrak{h} \otimes A \oplus \widetilde{H_N}$ with uniformly bounded dimensions. By Theorem \ref{T6.4}, $W^{\circ} \cong W_N$ as $\tau(H_N)(0)$-modules, where $\mathcal{Z}/\mathcal{K}$ acts trivially. Therefore the original module $V \cong V(W_N)$ which is the corresponding irreducible quotient with respect to the triangular decomposition given in (7.1). It is easy to see that $\mathcal{Z}/\mathcal{K}$ acts trivially using the fact that $V(W_N)$ is irreducible (use (2.5)(2)). We now provide an explicit realization of $V(W_N)$.

\bthm
Suppose that $V$ is a level zero irreducible integrable module for $\tau(H_N)$ having finite-dimensional weight spaces with respect to $\widetilde{\mathfrak{h}}$. Assume that $\mathfrak{g} \otimes A$ acts non-trivially on $V$. Then there exist a finite-dimensional irreducible module $V(\overline{\lambda})$ for $\mathfrak{g}$ and an irreducible finite-dimensional module $V_N$ over $\mathfrak{sp}_N$ such that $V \cong V(\overline{\lambda}) \otimes V_N \otimes A$. The action of $\tau(H_N)$ on the RHS is given by the following prescription.\\
$\widetilde{H_N}$ acts on the component $V_N \otimes A$ as described in Theorem \ref{T6.4}.\\
$X(r)(v \otimes w \otimes t^s) = (Xv) \otimes w \otimes t^{r+s} \ \forall \ X \in \mathfrak{g},  v \in V(\overline{\lambda}),  w \in V_N$ and $r,s \in \mathbb{Z}^N$.\\
$\mathcal{Z}/ \mathcal{K}$ acts trivially on $V$.  
\ethm

\begin{proof}
It is easy to check that the module at the RHS is irreducible over $\tau(H_N)$. We note that $\mathbb{C}v_{\lambda} \otimes V_N \otimes A$ is the top for $\tau(H_N)$ with respect to the triangular decomposition recorded in (7.1). Here $v_{\lambda}$ is a highest weight vector in $V(\overline{\lambda})$. We see that the top of $V$ is  $W^{\circ}$ and the top of RHS is $\mathbb{C}v_{\lambda} \otimes V_N \otimes A$, which are isomorphic as $\tau(H_N)(0)$-modules. Hence the corresponding irreducible quotients, which are $V$ and $V(\overline{\lambda}) \otimes V_N \otimes A$, are also isomorphic as $\tau(H_N)$-modules (see (2.14)). 
\end{proof}

\section{Contact Extended Affine Lie Algebra (KEALA)} \label{S8}
\noindent (8.1) In this section, we introduce yet another new class of extended affine Lie algebras. We call them by Contact extended affine Lie algebras (KEALAs) as the underlying derivation algebra resembles the Contact algebra (type $K$), KEALA is defined for odd number of variables. So let $M = 2m + 1, \ m \geqslant 1$.\\
(8.2) Let $\{e_1, e_2, \ldots, e_M \}$ be the standard basis and let $(\cdot, \cdot)$ be the standard bilinear form, i.e. $(e_i,e_j) = \delta_{ij}$. We shall treat $\mathbb{Z}^{2m}$ as a subgroup of $\mathbb{Z}^M$ with $(r,0) \in \mathbb{Z}^M$ for $r \in \mathbb{Z}^{2m}$.\\
(8.3) Let $r = (r_1, \ldots, r_m, r_{m+1}, \ldots, r_{2m}, r_M) \in \mathbb{Z}^M$. Define \\
$\overline{r} = (r_{m+1}, \ldots, r_{2m}, -r_1, -r_2, \ldots, -r_m) \in \mathbb{Z}^{2m} \subseteq \mathbb{Z}^M$ and also \\
$\underline{r} = \overline{r} + \sum_{i=1}^{2m} r_M e_i - \sum_{i=1}^{2m}r_i e_M$ \\
$= (r_{m+1} + r_M, \ldots, r_{2m} + r_M, -r_1 + r_M, \ldots, -r_m + r_M, - \sum_{i=1}^{2m} r_i) \in \mathbb{Z}^M$.\\
It is easy to check that $(\underline{r}, s) = (\overline{r}, s) + \sum_{i=1}^{2m}s_ir_M - \sum_{i=1}^{2m}r_is_M = - (r, \underline{s})$ so that $(\underline{r},r) = 0$.\\
(8.4) It is possible that $\underline{r} = 0$ even if $r \neq 0$. We first define a subgroup of $\mathbb{Z}^M$. Let $G = \{r \in \mathbb{Z}^M \ | \ r_i =  r_M, \ 1 \leqslant i \leqslant m \ \text{and} \ r_i = - r_M, \ m + 1 \leqslant i \leqslant 2m \}$. It is easy to see that $r \in G$ if and only if $\underline{r} = 0$.\\
(8.5) Define $D_M = \text{span} \{D(\underline{r},r) \ | \ r \in \mathbb{Z}^M \}$ and $\widetilde{D_M} = D_M \oplus D \subseteq DerA_M$. This definition is inspired by the work of Rudakov \cite{Ru}. It can be verified that for $r, s \in \mathbb{Z}^N$, $[D(\underline{r},r), D(\underline{s},s)] = (\underline{r},s) D(\underline{r} + \underline{s}, r + s)$.\\
We shall now define KEALA. Let \\
\begin{align*}
\mathcal{K}_M = \text{span} \big \{K(u,r) \ | \ u \in \mathbb{C}^M, \ 0 \neq r \in \mathbb{Z}^M, \ (u, \underline{r}) = 0 \big \}.
\end{align*}
It is trivial to check that $[D_M, K_M] \subseteq K_M$. Now consider 
\begin{align*}
\mathfrak{g} \otimes A_M \oplus \mathcal{Z} \oplus \widetilde{D_M}
\end{align*}
and note that $\mathcal{K}_M$ is an ideal in the above Lie algebra. Let us now define KEALA by setting
\begin{align*}
\tau(D_M) = \mathfrak{g} \otimes A_M \oplus \mathcal{Z}/\mathcal{K}_M \oplus \widetilde{D_M}.	
\end{align*} 

\bppsn
\
\begin{enumerate}
	
\item $\text{dim} (\mathcal{Z}/\mathcal{K}_M)_r = 0$ if $0 \neq r \in G$.
\item $\text{dim} (\mathcal{Z}/\mathcal{K}_M)_r = 1$ if $r \notin G$.
\item $\text{dim} (\mathcal{Z}/\mathcal{K}_M)_0 = M$.
\item $\text{dim} (D_M)_r = 0$ if $0 \neq r \in G$.
\item $\text{dim} (D_M)_r = 1$ if $r \notin G$.
\item  $\text{dim} (\widetilde{D_M})_0 = M$.   
\item $\tau(D_M)$ is an EALA.
\end{enumerate}
\eppsn

\begin{proof}
(1) and (3)-(6) are obvious. To prove (2), let $r \notin G$ so that $\underline{r} \neq 0$. Let $u, v \in \mathbb{C}^M$ such that $(u, \underline{r}) \neq 0 \neq (v, \underline{r})$. Consider
\begin{align*}
\bigg (\dfrac{u}{(u, \underline{r})} - \dfrac{v}{(v, \underline{r})}, \ \underline{r} \bigg ) = 0.
\end{align*}
This proves that  $K(u,r) = \lambda K(v,r)$ in $\mathcal{Z}/ \mathcal{K}_M$, where $\lambda$ is a non-zero scalar.\\
(*) We can take $K(\underline{s},s)$ as a basis for the vector space $(\mathcal{Z}/ \mathcal{K}_M)_s$ ($s \notin G$). \\ 
To prove (7), we only need to find an invariant, symmetric and non-degenerate bilinear form on $\tau(D_M)$. Note that $\widetilde{D_M} \subseteq S_M$ (see (3.4)). We shall now take the same form restricted to $\mathfrak{g} \otimes A_M \oplus \mathcal{Z} \oplus \widetilde{D_M}$. Note that $(D_M, K_M) =0$. Thus the form descends to $\tau(D_M)$. It can be easily verified that this form is invariant. For non-degeneracy, we shall also need to note that $(D(\underline{r},r), K(\underline{s},s)) = -\delta_{r+s,0}(\underline{s}, \underline{s})$, which is non-zero as $s \notin G$ (see (*)). 
\end{proof}

\brmk
One can define a twisted version of KEALA similar to the one defined in (3.8).
\ermk

\brmk
The EALA $\tau(S_N)$ defined in (3.8) in some sense is the largest possible EALA coming from an extension of a multiloop algebra. We also have a minimal EALA. Consider
\begin{align*}
\mathfrak{g} \otimes A \oplus \sum_{i=1}^{N} \mathbb{C}K_i \oplus \sum_{i=1}^{N} \mathbb{C}d_i,	
\end{align*}
with the Lie bracket defined as \\
(i) $[X(r), Y(s)] = [X,Y](r+s) + \delta_{r+s,0} <X,Y> \sum_{i =1}^{N}r_iK_i$, \\
(ii) $[K_i, d_j] = [K_i, K_j] = [d_i,d_j] = 0$, \\
(iii) $[d_i, X(r)] = r_iX(r)$, \\
(iv) $K_i$'s are central, \\
for all $X, Y \in \mathfrak{g}, \ r,s \in \mathbb{Z}^N, \ 1 \leqslant i,j \leqslant N$.

We also have an analogous notion for the twisted case. The irreducible integrable modules with finite-dimensional weight spaces for the minimal twisted EALA have been classified in \cite{S}, where each $K_i$ acts trivially.

If $N > 1$, then there does not exist any irreducible (not necessarily integrable) module with finite-dimensional weight spaces for this minimal twisted EALA, where at least one of the $K_i$'s acts non-trivially \cite{R}.
\ermk

\brmk
The irreducible integrable modules with finite-dimensional weight spaces over the twisted version of $\tau(S_N)$ have been classified in \cite{ESB} for the non-zero level case and in \cite{SB} for the level zero case. 
\ermk

\section{Appendix}
In this section, we prove a technical result which is used in Section \ref{S7} and Section \ref{S8}. We follow the approach of \cite{GL} where it is done for $N=2$. But there is a problem in the proof given in \cite{GL}. The authors of \cite{GL} have provided a revised version which we incorporate here, but our proof is for several variables. Our proof can be also used for $S_N$ instead of $H_N$.
(9.1) We maintain the same notations as in Section \ref{S3}. Recall that 
$H_N$ is a Hamiltonian Lie algebra and acts on $A = A_N$ via
\begin{align*}
h_r.t^s = (\overline{r},s)t^{r+s} \ \forall \ r, s \in \mathbb{Z}^N.
\end{align*}
Let $G = H_N \ltimes A$ be the semi-direct product of $H_N$ and $A$. Certainly $G = \bigoplus_{r \in \mathbb{Z}^N} G_r$ is a $\mathbb{Z}^N$-graded Lie algebra. Let $G^{\prime} = \bigoplus_{r \in \mathbb{Z}^N \setminus \{0 \}} G_r$, which is an ideal in $G$. Let $A^{\prime} = \bigoplus_{r \in \mathbb{Z}^N \setminus \{0 \}} \mathbb{C}t^r$.\\
Throughout this section, fix an irreducible $G$-module $W$. We assume that 
\begin{enumerate}
\item $W = \bigoplus_{r \in \mathbb{Z}^N} W_r$ and $G_r.W_s \subseteq W_{r+s}$.
\item $W$ is uniformly bounded in the sense that there exists $M \in \mathbb{N}$ such that dim $W_r \leqslant M \ \forall \ r \in \mathbb{Z}^N$.  
\end{enumerate} 

\bthm \label{T9.1}
Suppose that $A^{\prime}$ acts non-trivially on $W$ and let $t^0$ act as a non-zero scalar $c$ on $W$. Then there exist non-zero scalars $\lambda_{r,s}$ such that $t^rt^s = \lambda_{r+s}t^{r+s}$ on $W$ satisfying
\begin{enumerate}
\item $\lambda_{r,s} = \lambda \ \forall \ r, s \neq 0, \ r+s \neq 0$.
\item $\lambda_{r,-r} = \mu \ \forall \ r \neq 0$.
\item $\lambda_{0,r} = c$ and $\lambda^2 = \mu c$.
\end{enumerate}
\ethm 

We shall present the proof towards the end of this section. We first recall some standard results and prove some propositions.

The next two lemmas will follow from Lemma 3.1 and Lemma 3.2 in \cite{XLZ}. 

\blmma
Let $V = \bigoplus_{r \in \mathbb{Z}^N} V_r$ be an irreducible uniformly bounded $\mathbb{Z}^N$-graded module for $A$. Then dim $V_r \leqslant 1 \ \forall \ r \in \mathbb{Z}^N$.
\elmma

\blmma
Let $g \in U(A)$ such that $gv = 0$ for some $v \in W$. Then $g$ is locally nilpotent on $W$.
\elmma 

\bppsn \label{P9.4}
Either each $t^r, \ r \neq 0$, acts injectively on $W$ or every $t^r, \ r \neq 0$, acts locally nilpotently on $W$. 
\eppsn
\begin{proof}
We suppose that $t^r$ acts injectively on $W$ for some $r \neq 0$. Let $s \in \mathbb{Z}^N$ such that $(\overline{r},s) \neq 0$ and suppose that $t^s$ acts locally nilpotently on $W$. Then $t^s$ acts nilpotently on $W$. This follows using arguments similar to Claim 1 of \cite[Proposition 3.4]{GL}. Now using Claim 3 of \cite[Proposition 3.4]{GL} and the fact that $(\overline{r},s) \neq  0$, we get a contradiction. Hence it follows that $t^s$ is injective for all $s \in \mathbb{Z}^N$ satisfying $(\overline{r},s) \neq 0$. Now consider any $s \in \mathbb{Z}^N$ such that $(\overline{r},s) = 0$. Note that $(\overline{s} - r, r) \neq 0$. Hence $t^{s + \overline{r}}$ is injective. From above, it also follows that $t^{\overline{r}}$ is injective whenever $t^r$ is injective. Thus we conclude that $t^{\overline{s} - r}$ is injective. Now $(\overline{s}, \overline{s} - r) \neq 0$. Hence $t^s$ is also injective. This concludes the proof of this proposition.    
\end{proof}

\bppsn \label{P9.5}
If $t^r$ is locally nilpotent for all $r \neq 0$, then $A^{\prime}.W = 0$.
\begin{proof}
We closely follow the arguments given in \cite[Lemma 3.6]{GL}. Since $W$ is uniformly bounded, there exists some $M \in \mathbb{N}$ such that $(t^{-r}t^r)^MW = 0$ for all $r \neq 0$. Fix one such $r \neq 0$. The following can be deduced using arguments similar to the proof of \cite[Lemma 3.6]{GL}.
\begin{align*}
t^{s^1 - r} t^{s^2 - r} \ldots t^{s^{3M}-r} W = 0, \ \text{whenever} \ \sum_{i=1}^{3M} \epsilon_is^i \notin \mathbb{Q}r,  
\end{align*}
with the additional condition $(\overline{r},s^i) \neq 0$ ($s^i \in \mathbb{Z}^N, \ \epsilon_i \in \{0,1 \}$). \\
(3.7)-(3.10) of \cite{GL} can be established in our setup with the extra condition $(\overline{r},s^i) \neq 0$. Note that $\sum_{i=1}^{3M} \epsilon_i s_i$ are finitely many points in $\mathbb{Z}^N$. Now for any set of $n^i \in \mathbb{Z}^N \setminus \{0 \}, \ 1 \leqslant i \leqslant 3M$, choose $r \in \mathbb{Z}^N \setminus \{\sum_{i} \epsilon_im_i \}$ with $(\overline{r}, n^i) \neq 0$ (Such a choice of $r$ is possible. Let $M^i$ be the orthogonal complement of $n^i$. Choose $\overline{r} \in \mathbb{Z}^N \setminus \bigcup_{i}M^i$ such that $(\overline{r}, s^i) \neq 0$). Let $s^i = n^i + r$ and note that $(s^i, \overline{r}) \neq 0$. Now it is easy to verify that
\begin{align*}
\prod_{i=1}^{3M}t^{n^i}W = 0.
\end{align*} 
This holds good for any non-zero $\{n^i \ | \ 1 \leqslant i \leqslant 3M \}$. Finally consider the $G$-submodule $\{v \in W \ | \ t^rv = 0 \ \forall \ r \neq 0 \}$, which is non-zero and thus equal to $W$ from the irreducibility of $W$. 
\end{proof}
\eppsn

\noindent \textbf{Proof of Theorem \ref{T9.1}.} The proof runs parallel to \cite[Lemma 3.5]{GL}. Since $A^{\prime}$ acts non-trivially on $W$, from Proposition \ref{P9.4} and Proposition \ref{P9.5}, it follows that $t^r$ acts injectively on $W$ for all $r \neq 0$. This immediately implies that dim $W_r =$ dim $W_s = p$ (say) for all $r, s \in \mathbb{Z}^N$. As in \cite{GL}, define\\
(9.2)(1) $T_{r,s} = t^rt^s - \lambda_{r,s}t^{r+s} \ \forall \ r,s \in \mathbb{Z}^N$. \\
Then there exists a basis $\{w_1, \ldots, w_p \}$ of $W_0$ such that \begin{align*}
T_{r,s}(w_1 \ldots w_p) = t^{r+s}(w_1 \ldots w_p)(B_{r,s} - \lambda_{r,s}\text{Id}_p),
\end{align*}  
where $B_{r,s}$ is an upper triangular matrix with the diagonal entries equal to $\lambda_{r,s}$ and $\text{Id}_p$ is the identity matrix of order $p$.\\
Let us now record the following facts from \cite{GL}.\\
(9.2)(2) $T_{r,s}w_l \in \sum_{i=1}^{l-1} \mathbb{C}t^{r+s}w_i, \ r,s \in \mathbb{Z}^n, \ w_0 = 0$ and so $T_{r,s}^pw_p = 0$.\\
(9.2)(3) For any $r, s \in \mathbb{Z}^N$ and $u \in W_{-r-s} =  \sum_{i=1}^{p} \mathbb{C}t^{-r-s}w_i$, \\ $T_{r,s}u \subseteq \sum_{i=1}^{p-1} \mathbb{C}w_i$. Let $W_0^{\prime} = \sum_{i=0}^{p-1}\mathbb{C}w_i$.\\
We also note the following.\\
(9.2)(4) $t^{-r}t^r = T_{-r,r} + \lambda_{r, -r}t^0, \,\,\,\ t^{-r}t^rw_p \in \lambda_{r, -r}cw_p + W_0^{\prime}$.\\
Consider for some $l, k \in \mathbb{Z}^N$, $l = -(k+r+s)$. Then \\
$0 = (\prod t^l) (\prod h_k)T_{r,s}^p w_p$ (here we are taking product of $t^l$ exactly $p$ times and also the product of $h_k$ another $p$ times)
$= (\prod t^l) (\prod [h_k, T_{r,s}])w_p + T_{r,s}t^{-r-s}W_0 \equiv  (\prod t^l) (\prod [h_k, T_{r,s}])w_p (W_0^{\prime})$ (from (9.2)(3)) \\ 
$\equiv (\prod t^l) \prod \big ( (\overline{k}, r) t^{r+k}t^s + (\overline{k}, s)t^rt^{s+k} - (\overline{k}, r+s) \lambda_{r,s}t^{r+s+k} \big )w_p (W_0^{\prime})$ (follows from definitions) \\
$\equiv (\prod t^l) \prod \big ( (\overline{k}, s) (T_{r,k+s} + \lambda_{r, s+k}t^{r+s+k}) + (\overline{k}, r)(T_{r+k,s} + \lambda_{r+k,s}t^{r+s+k}) \\
- (\overline{k}, r+s) \lambda_{r,s}t^{r+s+k} \big )w_p (W_0^{\prime})$ (follows from definitions) \\
$\equiv (\prod t^l) \prod \big ( (\overline{k}, s) \lambda_{r,s+k} + (\overline{k}, r) \lambda_{s,r+k} - (\overline{k}, r+s) \lambda_{r,s} \big ) t^{r+s+k}w_p (W_0^{\prime})$.\\
Now using (9.2)(3) and (9.2)(4), we see that \\
$(\overline{k},s) \lambda_{r, k+s} + (\overline{k},r) \lambda_{s,r+s} - (\overline{k}, r+s) \lambda_{r,s} = 0$.\\
We shall now prove a series of lemmas to conclude Theorem \ref{T9.1}. We have \\
(9.3) For $l, r, s \in \mathbb{Z}^N \setminus \{0\}$, $(\overline{l}, s) \lambda_{r,s+l} + (\overline{l},r) \lambda_{s,r+l} - (\overline{l}, r+s) \lambda_{r,s} = 0$.

\blmma \label{L9.6}
For $(l, \overline{s}) \neq 0$, we have the following.
\begin{enumerate}
\item $\lambda_{s+l,s} = \lambda_{s,s}$.
\item $\lambda_{s,l} = \lambda_{s,s} = \lambda_{l,l}$.
\item $\lambda_{s,s} = \lambda_{s+l,s+l} = \lambda_{l,l}$.
\item $\lambda_{s+l,jl} = \lambda_{s,s} = \lambda_{l,l}, \ j \in \mathbb{Z} \setminus \{0 \}$.
\item $\lambda_{s, js+l} = \lambda_{s,s} = \lambda_{l,l}, \ j \in \mathbb{Z} \setminus \{0 \}$.
\end{enumerate}
\elmma

\begin{proof}
Take $r = js$ and $(\overline{l},s) \neq 0, \ j \neq 0$ in (9.3).\\
(9.3)(1) $\lambda_{js,s+l} + j \lambda_{js+l,s} = (j+1) \lambda_{js,s}$.\\
Now take $j=1$. Then $\lambda_{s,s+l} = \lambda_{s,s}$ and hence (1) follows.\\
Replace $l$ by $l-s$, then (2) follows. \\
From (2), we have $\lambda_{s+l,l} = \lambda_{s+l,s+l}$ and hence (3) follows.\\
To obtain (5), take $(j-1)s+l$ in place of $l$ in (1). \\ 
To see (4), we have $\lambda_{s+l,s} = \lambda_{s+l,s+l}$. Then $\lambda_{js+l,js} = \lambda_{js+l,js+l} \ \forall \ j \in \mathbb{Z}$. Again replace $l$ by $(1-j)s + l$ ($j \neq 1$) so that $\lambda_{l+s,js} = \lambda_{s,s} = \lambda_{l,l}$. 
\end{proof}   

\blmma \label{L9.7}
Let $r,s \in \mathbb{Z}^N \setminus \{0 \}, \ r \notin \mathbb{Q}s$ with $(r, \overline{s}) = 0$. Then there exists $k \in \mathbb{Z}^N$ such that $(r, \overline{k}) \neq 0$ and $(k, \overline{s}) = 0$.
\elmma

\begin{proof} First note that $(\overline{r}, \overline{s}) = (r,s)$.\\
\textbf{Case 1.} If $(r, s) = 0$, then simply take $\overline{k} = r$.\\
\textbf{Case 2.} If $(r,s) \neq 0$, then consider 
\begin{align*}
\overline{k} = \bigg (\dfrac{s}{(s,s)} - \dfrac{r}{(r,s)} \bigg )
\end{align*}
Now $(k, \overline{s}) = - (\overline{k},s) = 0$. Suppose that $(\overline{k},r) = 0$. Then it follows that $\dfrac{(s,r)}{(s,s)} = \dfrac{(r,r)}{(r,s)}$, which implies that $(s,r)^2 = (r,r)(s,s)$. It is a standard fact that $r \in \mathbb{Q}s$, which is a contradiction. This proves the lemma. 
\end{proof}

\blmma \label{L9.8}
Let $r \notin \mathbb{Q}s$. Then we have the following.
\begin{enumerate}
\item $\lambda_{r,r} = \lambda_{s,s}$.
\item $\lambda_{s+r,s} = \lambda_{s,s} = \lambda_{r,r}$.
\item $\lambda_{r,s} = \lambda_{s,s} = \lambda_{r,r}$.
\item $\lambda_{s+r,jr} = \lambda_{s,s} = \lambda_{r,r}$.
\item $\lambda_{s,js+r} = \lambda_{s,s} = \lambda_{r,r}$.
\end{enumerate}
\elmma

\begin{proof}
Suppose that $(s, \overline{r}) \neq 0$. Then the lemma follows from Lemma \ref{L9.6}. So let us assume that $(s, \overline{r}) = 0$. Then using Lemma \ref{L9.7}, pick $k \in \mathbb{Z}^N$ such that $(r, \overline{k}) \neq 0$ and $(k, \overline{s}) = 0$. Consequently $\lambda_{r,r} = \lambda_{r+k+s,r+k+s}$. But $(r+s, \overline{k}) \neq 0$ and henceforth (1) follows as $\lambda_{r,r} = \lambda_{r+s,r+s} = \lambda_{s,s}$.\\
To prove (2), we can assume that $(s, \overline{r}) = 0$. Take $k \in \mathbb{Z}^N$ such that $(r, \overline{k}) \neq 0$ and $(k, \overline{s}) = 0$. This implies that $\lambda_{r,r} = \lambda_{r+k+s,r}$. Now $(r+s, \overline{k}) \neq 0$ which gives $\lambda_{r+k+s,r} = \lambda_{r+s,r}$.\\
(4) and (5) follows from similar arguments as in Lemma \ref{L9.6}. 
\end{proof}

\blmma
\
\begin{enumerate}
\item $\lambda_{js,js} = \lambda_{s,s} \ \forall \ j \in \mathbb{Z} \setminus \{0 \}$.
\item $\lambda_{js,ps} = \lambda_{s,s}$ for $(j+p) \neq 0$.
\item $\lambda_{r,s} = \lambda \neq 0 \ \forall \ r+s \neq 0$ and $r, s \neq 0$. 
\end{enumerate}
\elmma

\begin{proof}
From Lemma \ref{L9.8}(1), we have $\lambda_{s,l} = \lambda_{s,s}$ for $l \notin \mathbb{Q}s$, which gives $\lambda_{js,l} = \lambda_{js,js}$. Applying Lemma \ref{L9.8}(4), we get $\lambda_{s,s} = \lambda_{js, s+l}$. Replacing $l$ by $l-s$, we obtain $\lambda_{s,s} = \lambda_{js,l} = \lambda_{js,js}$. To prove (2), considering $s$ as $ps$ ($p \neq 1$) in (9.3), we see that $p(\overline{l},s) \lambda_{r,ps+l} + (\overline{l},r) \lambda_{r+l,ps} = (\overline{l}, ps+r)$. Now choose $r=js$ and $(\overline{l},s) \neq 0$ to obtain that $p \lambda_{js,ps+l} + j \lambda_{js+l,ps} = (j+p) \lambda_{jr,ps}$. Consider $\lambda_{js,js+r} = \lambda_{js,js} = \lambda_{s,s}$. Take a suitable $r \in \mathbb{Z}^N$ so that
\begin{align*}
\lambda_{js, (j+p)s+r} = \lambda_{s,s} \ \text{where} \ j,p \ne 0, \ (j+p) \neq 0.
\end{align*}
Then $\lambda_{js,ps+r} = \lambda_{s,s}$ by choosing appropriate $r \in \mathbb{Z}^N$, which thus gives $(j+p) \lambda_{s,s} = (j+p) \lambda_{js,js}$ for $j+p \neq 0$. This implies that $\lambda_{s,s} = \lambda_{js,ps}$.   
\end{proof}

\blmma
\
\begin{enumerate}
\item $\lambda_{0,s} =c \ \forall \ s \in \mathbb{Z}^N$.
\item $\lambda_{r,-r} = \mu \ \forall \ r \neq 0$.
\item $c\mu = \lambda^2$. 
\end{enumerate}
\elmma
\begin{proof}
(1) Follows from the definition. \\
(2) Let $(r, \overline{s}) \neq 0$. Taking $l = -(s+r)$ in (9.3), we see that $\lambda_{r,-r} = \lambda_{s,-s}$. Suppose that $(r, \overline{s}) = 0$. Now one can choose $k \in \mathbb{Z}^N$ such that $(r, \overline{k}) \neq 0$ and $(s, \overline{k}) \neq 0$ (works for $N=2$ also). Then $\lambda_{r,-r} = \lambda_{k,-k} = \lambda_{s,-s}$.\\
(3) Consider $r,s, r+s \in \mathbb{Z}^N \setminus \{0 \}$. Then for all $w_1 \in W$,  
\begin{align*}
t^{-r}t^rt^st^{-s}w_1 = \lambda_{r,-r} \lambda_{s,-s}c^2w_1 = \mu^2c^2w_1, \\
t^{-r}t^{-s}t^rt^sw_1 = \lambda_{r,s} \lambda_{-r,-s} t^{r+s} t^{-(r+s)}w_1 = \lambda^2 \mu c w_1.
\end{align*}
This finally unravels that $\lambda^2 = \mu c$. 
\end{proof}

\end{document}